\documentclass[12pt,twoside]{article}%
\usepackage{enumerate}
\usepackage{enumitem}
\usepackage{amsmath,amssymb,amsfonts}
\usepackage{amsfonts}%
\usepackage{graphicx}
\usepackage{color}
\usepackage{pdfsync}
\providecommand{\U}[1]{\protect\rule{.1in}{.1in}}
\newtheorem{thm}{Theorem}[section]
\newtheorem{lem}[thm]{Lemma}

\newtheorem{prop}[thm]{Proposition}

\newtheorem{claim}[thm]{Claim}

\newenvironment{pf}[1][\bfseries Proof]{\noindent{#1.} }{\hfill \rule{0.5em}{0.5em}\\}
\newcommand{\fim}{\hfill\rule{2mm}{2mm}}
\numberwithin{equation}{section}
\begin{document}
\title{Multiplicity of solutions of  some quasilinear equations  in ${\mathbb{R}^{N}}$ with  variable exponents  
and concave-convex nonlinearities}

\author{Claudianor O. Alves\thanks{C.O. Alves was partially supported by CNPq/Brazil  301807/2013-2,  \hspace*{.7cm} e-mail:  coalves@dme.ufcg.edu.br}~~  Jos\'{e} L. P. Barreiro\thanks{J. L. P. Barreiro,~ e-mail: lindomberg@dme.ufcg.edu.br} \\
 Jos\'e V. A.  Goncalves\thanks{J.  V.  A.  Goncalves was partially supported by CNPq/PROCAD/UFG/UnB-Brazil, \hspace*{.5cm}   e-mail: goncalves.jva@gmail.com}}
\date{}
\maketitle

\begin{abstract}
In this paper, we prove multiplicity of solutions for a class of quasilinear problems in $ \mathbb{R}^{N} $ involving variable exponents and nonlinearities of concave-convex type. The main tools used are variational methods, more precisely, Ekeland's variational principle and  Nehari manifolds.
\end{abstract}



\section{Introduction}

In this paper, we consider the existence and multiplicity of solutions for the following class of quasilinear problems involving variable exponents
\begin{align}
\left\{
\begin{array}
[c]{rcl}%
-\Delta_{p(x)} u + \vert u \vert^{p(x) - 2} u & = & \lambda g(k^{-1}x)
\vert u \vert^{q(x) - 2}u + f(k^{-1} x) \vert u
\vert^{r(x) - 2}u~  \mbox{in}~ 
\mathbb{R}^{N},\\
u \in W^{1, p(x)}(\mathbb{R}^{N}), &  &
\end{array}
\right. \tag{$ P_{\lambda, k} $}\label{Plkm}%
\end{align}
where $\lambda$  and $ k $ are positive parameters with $ k \in \mathbb{N} $, the operator $ \Delta_{p(x)} u  = \mathrm{div}(|\nabla u|^{p(x) - 2} \nabla u) $  named $p(x)$-Laplacian,  is a natural extension of the $p$-Laplace operator, with $ p $ being a positive constant. We assume that $p,q, r : \mathbb{R}^{N} \to \mathbb{R} $ are positive Lipschitz
continuous functions, $ \mathbb{Z}^{N}$-periodic, that is, 
\begin{align}
p(x + z) = p(x), q(x + z) = q(x) \mbox{ and } r(x + z ) = r(x),~  x \in \mathbb{R}^{N}, \,\,   z \in \mathbb{Z}^{N}, \tag{$p_{1}$} \label{p1}
\end{align}
verifying
\begin{align}
1 <
q_{-} \leq q(x) \leq q_{+} < p_{-} \leq p(x) \leq p_{+} < r_{-} \leq r \ll p^{*} \tag{$p_{2}$} \label{p2},~ \mbox{a.e. in}~  \mathbb{R}^{N}, 
\end{align}
where $ p_{+} = \mathrm{ess \;sup}_{x \in \mathbb{R}^{N}} p(x)$,  $ p_{-} = \mathrm{ess \;inf}_{x \in \mathbb{R}^{N}} p(x)$ and

\begin{equation}
p^{*}(x)  = \left\{
\begin{array}{cl}
Np(x)/(N - p(x))& \mbox{ if }  p(x) < N,\\
+ \infty & \mbox{ if } p(x) \geq N.
\end{array}
\right. \tag{$P$}\label{probP}
\end{equation}
Hereafter, the notation $u \ll v$ means that $ \displaystyle \inf_{x \in \mathbb{R}^{N}} (v(x) - u(x)) > 0$.

Furthermore, we assume the  condition:
\begin{enumerate}[label={(H)}]
\setcounter{enumi}{0}
\item\label{H4} \hspace{2 cm} $\displaystyle \frac{q_{+}}{p_{-}} < \frac{(r_{+} - q_{+})}{(r_{+} - p_{-})} \frac{(r_{-} - p_{+})}{(r_{-} - q_{-})}$.
\end{enumerate}
Here, we would like to point out that this condition  is equivalent to $0<q<p$ for the case where the exponent is constant. This technical condition will be needed, especially  in the proof of Lemma \ref{Jlk-Mlk-menos}.

\vspace{0.5 cm}

Regarding the functions $f$ and $g $, we assume the following conditions:

\begin{enumerate}[label={($g_\arabic{*}$})]
\setcounter{enumi}{0}
\item\label{H1} $ g : \mathbb{R}^{N} \to \mathbb{R} $ is a nonnegative measurable function with $ g \in L^{\Theta(x)}(\mathbb{R}^{N}) $ where $ \Theta(x) = \frac{r(x)}{r(x) - q(x)}$,
\end{enumerate}
\begin{enumerate}[label={($f_\arabic{*}$})]
\setcounter{enumi}{0}
\item\label{H2} $ f : \mathbb{R}^{N} \to \mathbb{R} $ is  a positive continuous function  such that
\[
\lim_{\vert x \vert\rightarrow\infty} f(x) = f_{\infty}
\]
and $0 < f_{\infty} < f(x) $ for all $x \in\mathbb{R}^{N} $,
\item\label{H3}
there exist $\ell$ points $a_{1}, a_{2}, \cdots,a_{\ell} $ in
$\mathbb{Z}^{N} $ with $ a _{1} = 0 $, such that
\[
1 = f(a_{i}) = \max_{\mathbb{R}^{N}}f(x), \text{ for }1 \leq i \leq\ell .
\]
\end{enumerate}

Problems with variable exponents  appear in various applications. The reader is referred to  R\r{u}\v{z}i\v{c}ka \cite{Ru} and Krist\'aly, Radulescu \& Varga in \cite{KRV} for  several questions in mathematical physics where such class of problems appear. In recent years,these  problems have attracted an increasing attention.  We would like to mention \cite{alves08, AlvesFerreira2, AlvesSouto,AlvesShibo, chabrowki, FanHan,fu11, MR}, as well as  the survey papers \cite{AS,DHN,S} for the advances and references in this field.

The problem $(P_{\lambda,k})$ has been considered in the literature for the case where the exponents are constants, see for example, Adachi \& Tanaka \cite{AT}, Cao \& Noussair \cite{CN}, Cao \& Zhou \cite{Cao}, Hirano \cite{H1}, Hirano \& Shioji \cite{HS}, Hu \& Tang \cite{HuTang}, Jeanjean \cite{jeanjean}, Lin \cite{Lin12}, Hsu, Lin \& Hu \cite{Hsu1}, Tarantello \cite{T}, Wu \cite{Wu1, Wu2} and their references.

\medskip

In Cao \& Noussair \cite{CN}, the authors have studied the existence and multiplicity of positive and nodal solutions for the following problem
$$
\left\{
\begin{array}{l}
-\Delta u + u  =  f(\epsilon x) \vert u \vert^{r - 2}u  \,\,\, \mbox{in} \,\,\, \mathbb{R}^{N}\\
\mbox{}\\
u \in H^{1,2}(\mathbb{R}^{N}),
\end{array}
\right. \eqno{(P_{1})}
$$
where $\epsilon$ is a positive real parameter, $r \in (2,2^{*})$  and $f$ verifies conditions \ref{H2}-\ref{H3}. By using variational methods, the authors showed the existence of at least $\ell$ positive solutions and $\ell$ nodal solutions if $\epsilon$ is small enough. Later on, Wu in \cite{Wu1} considered the perturbed problem
$$
\left\{
\begin{array}{l}
-\Delta u + u  =  f(\epsilon x) \vert u \vert^{r - 2}u + \lambda g(\epsilon x)|u|^{q-2}u \,\,\, \mbox{in} \,\,\, \mathbb{R}^{N}\\
\mbox{}\\
u \in H^{1,2}(\mathbb{R}^{N}),
\end{array}
\right. \eqno{(P_{2})}
$$
where $\lambda$ is a positive parameter and $q \in (0,1)$. In \cite{Wu1}, the authors showed the existence of at least $\ell$ positive solutions for $(P_2)$ when  $\epsilon$ and $\lambda$ are small enough.

\medskip

In Hsu, Lin \& Hu \cite{Hsu1}, the authors have considered the following class of quasilinear problems
$$
\left\{
\begin{array}{l}
-\Delta_p u + |u|^{p-2}u  =  f(\epsilon x) \vert u \vert^{r - 2}u + \lambda g(\epsilon x) \,\,\, \mbox{in} \,\,\, \mathbb{R}^{N}\\
\mbox{}\\
u \in W^{1,p}(\mathbb{R}^{N})
\end{array}
\right. \eqno{(P_{3})}
$$
with $N \geq 3$ and $2 \leq p < N$. In that paper, the authors have proved the same type of results found in \cite{CN} and \cite{Wu1}.

\medskip

Motivated by results proved in \cite{CN}, \cite{Hsu1} and \cite{Wu1}, we intend in the present paper to prove the existence of multiple solutions for problem (\ref{Plkm}), by using the same type of approach explored in those papers. However, once that we are working with variable exponents, some estimates that hold for the constant case are not immediate for the variable case, and so, a careful analysis is necessary to get some estimates. More precisely, when the exponents are constant each term in the nonlinearity is homogeneous, which is very good to get some estimates involving the energy functional, however if the exponents are not constant we loose this property. Here,  this difficulty is overcome by using Lemmas \ref{lemAbst} and \ref{Max-Nehari}.  We added further explanations  immediately before the statement of each of these lemmas.

\bigskip
Our main result is the following

\begin{thm}\label{T1}
Assume that (\ref{p1})--(\ref{p2}), \ref{H1}, \ref{H2}--\ref{H3} and \ref{H4} are satisfied. Then,
there are positive numbers $ k_{*} $ and $ \Lambda_{*}=\Lambda(k_{*}) $, such that problem (\ref{Plkm}) admits at least $ \ell + 1 $ solutions for $ 0 < \lambda < \Lambda_{*} $ and $ k > k_{*}$.
\end{thm}

\vspace{0.5 cm}
\noindent\textbf{Notation:} The following notations will be used in the present work:
\begin{itemize}
  \item $ C $ and $ c_{i} $ denote generic positive constants, which may vary from line to line.
  \item  We denote by $\int u $ the integral $\int_{\mathbb{R}^{N}}udx$, for any measurable function $u$.
  \item $ B_{R}(z) $ denotes the open ball with center at $ z $ and radius $ R $ in $\mathbb{R}^{N}$.
\end{itemize}

\section{Preliminaries on Lebesgue and Sobolev spaces with variable exponents  in $\mathbb{R}^{N}$}

\label{expoents_variaveis} In this section, we recall the definitions and some
results involving the spaces $L^{h(x)}(\mathbb{R}^{N}) $ and $W^{1,h(x)}%
(\mathbb{R}^{N}) $.  We refer to \cite{peter,Fan2001a, Fan2001b, kovacik91} for the fundamental properties of these spaces.

Hereafter, let us denote by $L_{+}^{\infty}(\mathbb{R}^{N})$ the set
\[
L_{+}^{\infty}(\mathbb{R}^{N}) = \left\{  u \in L^{\infty}(\mathbb{R}^{N}) : \mbox{ess}\inf_{x
\in\mathbb{R}^{N}}u \geq1\right\},
\]
and we will assume that $h \in L_{+}^{\infty}(\mathbb{R}^{N})$.

The variable exponent Lebesgue space $L^{h(x)}(\mathbb{R}^{N}) $ is defined by
\[
L^{h(x)}(\mathbb{R}^{N}) = \left\{  u \colon\mathbb{R}^{N} \to \mathbb{R} \text{ is measurable } \colon
\, \, \,\, \int \vert u(x)\vert^{h(x)} < + \infty\right\},
\]
and its usual norm is
\[
\Vert u \Vert_{h(x)} = \inf\left\{  t > 0 : \int
\left\vert \frac{u(x)}{t}\right\vert ^{h(x)} \leq1\right\} .
\]

On the space $L^{h(x)}(\mathbb{R}^{N})$, we consider the \textit{modular function} $\rho:
L^{h(x)}(\mathbb{R}^{N}) \to\mathbb{R}$ given by
\[
\rho(u) = \int |u(x)|^{h(x)} .
\]

In what follows, let us  denote by $h_{-}$ and $h_{+}$ the following real
numbers
\[
h_{-} = \mbox{ess}\inf_{x \in \mathbb{R}^{N}}{h(x)} \,\,\,\,\, and \,\,\,\,\,  h_{+} =  \mbox{ess}\sup_{x \in \mathbb{R}^{N}}{h(x)}.
\]

\begin{prop}
\label{modular} Let $u \in L^{h(x)}(\mathbb{R}^{N}) $ and $\{u_{n}\}_{n \in\mathbb{N}}
\subset L^{h(x)}(\mathbb{R}^{N}) $. Then,

\begin{enumerate}
\item If $u \neq0 $, $\Vert u \Vert_{h(x)} = a
\Leftrightarrow\rho\left(  \frac{u}{a}\right)  = 1$.

\item $\Vert u \Vert_{h(x)} < 1 \quad(=1; > 1) \Leftrightarrow
\rho(u) < 1 (= 1; > 1) $;

\item $\Vert u \Vert_{h(x)} > 1 \Rightarrow\Vert u \Vert_{h(x)}^{h_{-}} \leq\rho(u) \leq\Vert u \Vert_{h(x)}^{h_{+}} $.

\item $\Vert u \Vert_{h(x)} < 1 \Rightarrow\Vert u \Vert_{h(x)}^{h_{+}} \leq\rho(u) \leq\Vert u \Vert_{h(x)}^{h_{-}} $.

\item $\displaystyle \lim_{n \to+\infty} \Vert u_{n} \Vert_{h(x)}
= 0 \Leftrightarrow\lim_{n \to\infty}\rho(u_{n}) = 0 .$

\item $\displaystyle \lim_{n \to+\infty} \Vert u_{n} \Vert_{h(x)}
= + \infty\Leftrightarrow\lim_{n \to\infty}\rho(u_{n}) = + \infty$.
\end{enumerate}
\end{prop}

As usual, we denote by $h^{\prime}(x) = \frac{h(x)}{h(x) - 1} $ the conjugate exponent function
of $h(x) $, and define
\begin{align*}
h^{*}(x) = \left\{
\begin{array}
[c]{lcc}%
\dfrac{Nh(x)}{N - h(x)} & \text{if} & h(x) < N\\
+ \infty & \text{if} & h(x) \geq N.
\end{array}
\right.
\end{align*}

We have the following H\"older inequality for Lebesgue spaces with variable exponents.
\begin{prop}
[H\"{o}lder-type Inequality]Let $u \in L^{h(x)}(\mathbb{R}^{N})$ and $v \in L^{h^{\prime
}(x)}(\mathbb{R}^{N}) $. Then, $uv \in L^{1}(\mathbb{R}^{N})$ and
\begin{align*}
\int \vert u(x)v(x)\vert  \leq\left(  \frac{1}{h_{-}} + \frac
{1}{h^{\prime}_{-}} \right)  \Vert u \Vert_{h(x)} \Vert v
\Vert_{h^{\prime}(x)}.
\end{align*}

\end{prop}

\begin{lem}
\label{modular:interp} Let $h,b \in L_{+}^{\infty}(\mathbb{R}^{N})$ with $h(x) \leq
b(x)$ a.e. in $\mathbb{R}^{N}$ and $u \in L^{b(x)}(\mathbb{R}^{N}) $. Then, $|u|^{h(x)}
\in L^{\frac{b(x)}{h(x)}}(\mathbb{R}^{N})$,
\begin{align*}
\Vert\vert u \vert^{h(x)}\Vert_{\frac{b(x)}{h(x)}} \leq\max\left\{
\Vert u \Vert^{h_{+}}_{b(x)}, \Vert u \Vert^{h_{-}}_{b(x)}\right\},
\end{align*}
and further
\begin{align*}
\Vert\vert u \vert^{h(x)} \Vert_{\frac{b(x)}{h(x)}} \leq\Vert u
\Vert^{h_{+}}_{b(x)} + \Vert u \Vert^{h_{-}}_{b(x)}.
\end{align*}
\end{lem}

The next three results are important tools to study the properties of some energy functionals, and their proofs can be found in \cite{AlvesFerreira2}.


%

\begin{prop} [Brezis-Lieb's lemma, first version] \label{Brezis-Lieb-1}
   Let $ \{ \eta_n \} \subset L^{ h(x) } ( \mathbb R^N, \mathbb R^m ) $ with $ m \in \mathbb{N} $ verifying

   \begin{enumerate}
      \item [\emph{(i)}] $ \eta_n(x) \to \eta(x), \ \text{a.e. in} \ \mathbb R^N $;
      \item [\emph{(ii)}] $ \displaystyle \sup_{n \in \mathbb N } | \eta_n |_{ L^{ h(x) } ( \mathbb R^N, \mathbb R^m ) } < \infty $. \\
   \end{enumerate}
   Then, $ \eta \in L^{ h(x) } ( \mathbb R^N, \mathbb R^m ) $ and
   \begin{equation}
      \int \left( \left| \eta_n \right|^{ h(x) } - \left| \eta_n - \eta \right|^{ h(x) } - \left| \eta \right|^{ h(x) } \right) \,= o_n(1).
   \end{equation}
\end{prop}

\begin{prop} [Brezis-Lieb's lemma, second version] \label{Brezis-Lieb-2}
   Let $\{ \eta_n \} \subset L^{ h(x) } ( \mathbb R^N, \mathbb R^m ) $ verifying

   \begin{enumerate}
      \item [\emph{(i)}] $ \eta_n(x) \to \eta(x), \ \text{a.e. in} \ \mathbb R^N $;
      \item [\emph{(ii)}] $ \displaystyle \sup_{n \in \mathbb N } | \eta_n |_{ L^{ h(x) } ( \mathbb R^N, \mathbb R^m ) } < \infty $. \\
   \end{enumerate}
   Then
   \begin{equation}
      \eta_n \rightharpoonup \eta \ \text{in} \ L^{ h(x) } ( \mathbb R^N, \mathbb R^m ).
   \end{equation}
\end{prop}

The next proposition is a Brezis-Lieb type result.

\begin{prop} [Brezis-Lieb lemma, third version]
   Let   $ \{ \eta_n \}\subset L^{ h(x) }  ( \mathbb R^N, \mathbb R^m ) $  such that

   \begin{enumerate} \label{Brezis-Lieb-3}
      \item [\emph{(i)}] $ \eta_n(x) \to \eta(x), \ \text{a.e. in} \ \mathbb R^N $;
      \item [\emph{(ii)}] $ \displaystyle \sup_{n \in \mathbb N } | \eta_n |_{ L^{ h(x) } ( \mathbb R^N, \mathbb R^m ) } < \infty $. \\
   \end{enumerate}
Then
   \begin{equation}
      \int \left| \left| \eta_n \right|^{ h(x)-2 } \eta_n - \left| \eta_n - \eta \right|^{ h(x)-2 } \left( \eta_n - \eta \right) - \left| \eta \right|^{ h(x)-2 } \eta \right|^{ h'(x) } \,  = o_n(1).
   \end{equation}
 \end{prop}

The variable exponent Sobolev space $W^{1,h(x)}(\mathbb{R}^{N}) $ is defined by
\begin{align*}
W^{1,h(x)}(\mathbb{R}^{N}) = \left\{  u \in W^{1,1}_{loc}(\mathbb{R}^{N}) : u \in
L^{h(x)}(\mathbb{R}^{N}) \quad\text{ and } \quad| \nabla u | \in L^{h(x)}(\mathbb{R}^{N})
\right\}.
\end{align*}
The corresponding norm for this space is
\begin{align*}
\Vert u \Vert_{W^{1,h(x)}(\mathbb{R}^{N})} = \Vert u \Vert_{h(x)} +
\Vert\nabla u \Vert_{h(x)}.%
\end{align*}
The spaces $L^{h(x)}(\mathbb{R}^{N})$ and $W^{1,h(x)}(\mathbb{R}^{N})$ are separable and reflexive Banach spaces when $h_{-} >1$.

On the space $W^{1,h(x)}(\mathbb{R}^{N})$, we consider the {\it modular function}  \linebreak $\rho_{1}: W^{1,h(x)}(\mathbb{R}^{N}) \to \mathbb{R}$ given by
\[
\rho_{1}(u) = \int \left(  |\nabla u(x)| ^{h(x)} + |u(x)| ^{h(x)}\right).
\]
If, we define
\begin{align}\label{norm}
\Vert u \Vert = \inf\left\{ t > 0 : \int \frac{(|\nabla u |^{h(x)}+| u |^{h(x)})}{t^{h(x)} }  \leq 1\right\},
\end{align}
then $ \Vert \cdot \Vert_{W^{1,h(x)}(\mathbb{R}^{N})} $ and $ \Vert \cdot \Vert $ are equivalent norms on $ W^{1,h(x)}(\mathbb{R}^{N}) $.
\begin{prop}\label{modular2}
Let $ u \in W^{1,h(x)}(\mathbb{R}^{N}) $ and $ \{u_{n}\} \subset W^{1,h(x)}(\mathbb{R}^{N}) $. Then, the same conclusion of Proposition \ref{modular} occurs considering  $\|\cdot\|$ and $\rho_{1}$.

\end{prop}
%
%

%

The next result is a Sobolev embedding Theorem for variable exponent, whose proof can be found in
\cite{peter} and \cite{Fan2001a}.

\begin{thm} [Sobolev embedding] \label{timersao}
Let $ h : \mathbb{R}^{N} \to \mathbb{R} $ be  Lipschitz continuous with  $1 < h_{-} \leq h_{+} < N$, and consider $ s \in L^{\infty}_{+}(\mathbb{R}^{N}) $ satisfying  $ h(x) \leq s(x) \leq h^{*}(x) $ a.e. in  $ \mathbb{R}^{N} $. Then there is the  continuous embedding
$$
W^{1,h(x)}(\mathbb{R}^{N}) \hookrightarrow  L^{s(x)}(\mathbb{R}^{N}).
$$
\end{thm}

\section{Technical lemmas}
For convenience, in all this paper, we define the following functions
\begin{align*}
g_{k}(x) =  g(k^{-1} x) \quad \mbox{ and } \quad f_{k}(x) =  f(k^{-1} x) \quad \mbox{ for all } x \in \mathbb{R}^{N}.
\end{align*}

Associated with the problem~(\ref{Plkm}), we have the energy functional
\linebreak$J_{\lambda, k}:W^{1,p(x)}(\mathbb{R}^{N})\rightarrow\mathbb{R}$ defined
by
\begin{align*}
J_{\lambda, k}(u) &  =\int\frac{1}{p(x)}\left(  |\nabla
u|^{p(x)}+|u|^{p(x)}\right)  -\lambda\int\frac{g_{k}(x)}{q(x)}|u|^{q(x)}
 - \int  \frac{f_{k}(x)}{r(x)}|u|^{r(x)} .
\end{align*}
A direct computation gives $J_{\lambda, k}\in C^{1}\left(  W^{1,p(x)}%
(\mathbb{R}^{N}),\mathbb{R}\right)  $ with
\begin{align*}
J'_{\lambda, k}(u) v  & =\int\left(  |\nabla
u|^{p(x)-2}\nabla u\nabla v+|u|^{p(x)-2}uv\right)  -\lambda\int g_{k}(x)|u|^{q(x)-2}uv\\
& \quad -\int f_{k}(x) |u|^{r(x)-2}uv  ,
\end{align*}
for each  $u,v\in W^{1,p(x)}(\mathbb{R}^{N})$. Therefore, the critical points of $ J_{\lambda, k} $ are precisely the (weak) solutions of  (\ref{Plkm}).

Since $ J_{\lambda, k} $ is not bounded from below on $ W^{1,p(x)}(\mathbb{R}^{N}) $, we will work on the \emph{Nehari manifold} $\mathcal{M}_{\lambda, k}$ associated  with $ J_{\lambda, k} $, given by
$$
\mathcal{M}_{\lambda, k} = \left\{ u \in W^{1,p(x)}(\mathbb{R}^{N})\setminus\{ 0\}: J'_{\lambda, k}(u) u = 0 \right\}.
$$
In what follows, we denoted by $c_{\lambda, k}$ the  real number
\[
c_{\lambda, k} = \inf_{u \in \mathcal{M}_{\lambda, k}} J_{\lambda, k}(u).
\]
Using well-known arguments, it is easy to prove that $c_{\lambda, k} $ is the mountain pass level of  $J_{\lambda, k}$.

For $ f \equiv 1 $ and $ \lambda = 0 $, we consider the problem
\begin{align}
\left\{
\begin{array}
[c]{rcl}%
-\Delta_{p(x)} u + \vert u \vert^{p(x) - 2} u & = & \vert u
\vert^{r(x) - 2}u , \quad
\mathbb{R}^{N}\\
\mbox{}\\
u \in W^{1, p(x)}(\mathbb{R}^{N}). &  &
\end{array}
\right. \tag{$ P_{\infty} $}\label{Poo}%
\end{align}
Associated with the problem~(\ref{Poo}), we have the energy functional \linebreak $ J_{\infty}:W^{1, p(x)}(\mathbb{R}^{N}) \to \mathbb{R}$ given by
\begin{align*}
J_{\infty}(u)   =\int \frac{1}{p(x)}\left(  |\nabla
u|^{p(x)}+|u|^{p(x)}\right)    -  \int \frac{1
}{r(x)}|u|^{r(x)},
\end{align*}
the mountain pass level
$$
c_{\infty} = \inf_{u \in \mathcal{M}_{\infty}} J_{\infty}(u),
$$
and the Nehari manifold
$$
\mathcal{M}_{\infty} = \left\{ u \in W^{1,p(x)}(\mathbb{R}^{N})\setminus\{ 0\}: J'_{\infty}(u)u = 0 \right\}.
$$

For $ f \equiv f_{\infty} $ and $ \lambda = 0 $, we fix the problem
\begin{align}
\left\{
\begin{array}
[c]{rcl}%
-\Delta_{p(x)} u + |u|^{p(x) - 2} u & = & f_{\infty} |u|^{r(x) - 2}u, \quad
\mathbb{R}^{N}\\
\mbox{}\\
u \in W^{1, p(x)}(\mathbb{R}^{N}), &  &
\end{array}
\right. \tag{$ P_{f_{\infty}} $}\label{Pf00}%
\end{align}
and as above, we denote by $ J_{f_\infty}, c_{f_{\infty}}$ and $\mathcal{M}_{f_{\infty}}$ the energy functional, the mountain pass level, and Nehari manifold associated with $ (P_{f_{\infty}}) $ respectively.

\medskip

Hereafter, let us fix $ K > 1 $ such that
\begin{equation}\label{K}
\Vert v \Vert_{r(x)} \leq K \Vert v \Vert \mbox{ for any }  v \in  W^{1, p(x)}(\mathbb{R}^{N}),
\end{equation}
which exists by Theorem \ref{timersao}.

\vspace{0.5 cm}

The next lemma is a technical result, which will be used in Section 5.

\begin{lem}[Local property] \label{local}
For each $ k \in \mathbb{N} $, there are positive constants $ \lambda^{*} = \lambda^{*}(k),\beta$ and $ \sigma $ independent of $k$, such that $ J_{\lambda, k}(u) \geq \beta > 0 $ for all $ \lambda \in (0, \lambda^{*}) $ with  $ \Vert u \Vert = \sigma $.
\end{lem}

\begin{pf}
Combining the definition of $ J_{\lambda, k} $ with H\"{o}lder's inequality, Sobolev embedding and Proposition~\ref{modular:interp}, we derive
\begin{align*}
J_{\lambda, k}(u) \geq & \frac{1}{p_{+}} \int \left(  |\nabla u|^{p(x)}+|u|^{p(x)}\right)   - 2\frac{\lambda}{q_{-}} \Vert g_{k} \Vert_{\Theta(x)} K^{q_{+}} \max\{ \Vert u \Vert^{q_{-}}, \Vert u \Vert^{q_{+}}\}  \\
& \quad - \int \frac{1}{r_{-}} |u|^{r(x)}.
\end{align*}

\noindent If $ \Vert u \Vert < 1 $,  by Proposition~\ref{modular2} and
Theorem \ref{timersao},
\begin{align*}
J_{\lambda, k}(u) \geq & \frac{1}{p_{+}} \Vert u \Vert^{p_{+}} - 2\frac{\lambda}{q_{-}} \Vert g_{k} \Vert_{\Theta(x)} K^{q_{+}} \Vert u \Vert^{q_{-}} -  \frac{K^{r_{+}}}{r_{-}} \Vert u \Vert^{r_{-}}.
\end{align*}
Since $ p_{+} < r_{-} $, by fixing $ \sigma$ small enough such that
\[
\frac{1}{p_{+}} \sigma^{p_{+}} - \frac{K^{r_{+}}}{r_{-}} \sigma^{r_{-}} \geq \frac{1}{2p_{+}} \sigma^{p_{+}},
\]
we obtain
\begin{align*}
J_{\lambda, k}(u) \geq & \frac{1}{2 p_{+}} \sigma^{p_{+}} -2 \frac{\lambda}{q_{-}} \Vert g_{k} \Vert_{\Theta(x)} K^{q_{+}} \sigma^{q_{-}},
\end{align*}
for $ \Vert u \Vert = \sigma $.

Now, fix $ \lambda^{*} = \lambda^{*}(k) > 0 $ satisfying
\begin{align}
\lambda^{*} \Vert g_{k} \Vert_{\Theta(x)} < \frac{q_{-}}{8p_{+}K^{q_{+}}} \sigma^{p_{+} - q_{-}}. \label{local_prop}
\end{align}
Then, if $ 0 < \lambda < \lambda^{*} $,
\[
J_{\lambda, k}(u) > \frac{1}{2p_{+}} \sigma^{p_{+}} - 2\frac{\lambda^{*}}{q_{-}} \Vert g_{k} \Vert_{\Theta(x)} K^{q_{+}} \sigma^{q_{-}} >  \frac{1}{4p_{+}} \sigma^{p_{+}} = \beta >  0 \text{ on } \partial B_{\sigma}(0),
\]
proving the result.
\end{pf}

The next result concerns with the behavior of $ J_{\lambda, k} $ on $ \mathcal{M}_{\lambda, k} $.
\begin{lem}
The energy functional $ J_{\lambda, k} $ is coercive and bounded from below on $ \mathcal{M}_{\lambda, k} $.
\end{lem}
\begin{pf}
For $ u \in \mathcal{M}_{\lambda, k} $, we have $ J'_{\lambda, k}(u)  u = 0 $. Therefore,
\begin{align*}
\int f_{k}(x) |u|^{r(x)}   = \int \left(  |\nabla
u|^{p(x)}+|u|^{p(x)}\right)   - \lambda \int g_{k}(x) |u|^{q(x)},
\end{align*}
loading to
\begin{align*}
J_{\lambda, k}(u) &\geq  \frac{1}{p_{+}} \int \left(  |\nabla u|^{p(x)}+|u|^{p(x)}\right)   - \frac{\lambda}{q_{-}} \int g_{k}(x) |u|^{q(x)} -  \frac{1}{r_{-}}\int f_{k}(x)  |u|^{r(x)} \\
&=  \left( \frac{1}{p_{+}} - \frac{1}{r_{-}}\right) \int \left(  |\nabla u|^{p(x)}+|u|^{p(x)}\right)    - \lambda \left( \frac{1}{q_{-}} - \frac{1}{r_{-}}\right) \int g_{k}(x) |u|^{r(x)}.
\end{align*}
If $ \Vert u \Vert  > 1 $, the  Propositions \ref{modular:interp} and \ref{modular2} together with  H\"{o}lder's inequality and Theorem \ref{timersao} give
\begin{align*}
J_{\lambda, k}(u)  \geq &\left( \frac{1}{p_{+}} - \frac{1}{r_{-}}\right) \Vert u \Vert^{p_{-}}   - \left( \frac{1}{q_{-}} - \frac{1}{r_{-}}\right) 2 \lambda K^{q^{+}} \Vert g_{k} \Vert_{\Theta(x)} \Vert u \Vert^{q_{+}} \\
 = & \Vert u \Vert^{q_{+}} \left\{ \left( \frac{1}{p_{+}} - \frac{1}{r_{-}}\right) \Vert u \Vert^{p_{-} - q_{+}} - 2 \lambda \left( \frac{1}{q_{-}} - \frac{1}{r_{-}}\right) K^{q^{+}} \Vert g_{k} \Vert_{\Theta(x)}  \right\}.
\end{align*}
Since $ q_{+} < p_{-} $, the last inequality implies that $ J_{\lambda, k}$ is coercive and bounded from below on  $ \mathcal{M}_{\lambda, k} $.
\end{pf}

From now on, let
\[ E_{\lambda, k}(v) = J'_{\lambda, k}(v) v \quad \mbox{ for any } v \in W^{1, p(x)}(\mathbb{R}^{N}) .
\]
Employing the functional $E_{\lambda,k}$, we split $ \mathcal{M}_{\lambda, k} $ into three parts:
\begin{align*}
\mathcal{M}_{\lambda, k}^{+} = \{ v \in \mathcal{M}_{\lambda, k} : E'_{\lambda, k}(v)v > 0 \}, \\[0.5cm]
\mathcal{M}_{\lambda, k}^{0} = \{ v \in \mathcal{M}_{\lambda, k} : E'_{\lambda, k}(v)v = 0 \}, \end{align*}
and
$$
\mathcal{M}_{\lambda, k}^{-} = \{ v \in \mathcal{M}_{\lambda, k} : E'_{\lambda, k}(v)v < 0 \}.
$$

In the next lemma, we prove that the critical points of $ J_{\lambda, k} $ restrict to $ \mathcal{M}_{\lambda, k} $ which do not belong $ \mathcal{M}^{0}_{\lambda, k} $ are in fact critical points of $ J_{\lambda, k} $ on $ W^{1,p(x)}(\mathbb{R}^{N})$.
\begin{lem}
If $ u_{0} \in \mathcal{M}_{\lambda, k} $ is a critical point of $ J_{\lambda, k} $ restricted to  $ {\mathcal{M}_{\lambda, k}} $ and $ u_{0} \not\in \mathcal{M}_{\lambda, k}^{0} $, then $ u_{0} $ is a critical point of $ J_{\lambda, k}$.
\end{lem}

\begin{pf}
By Lagrange multiplier theorem, there is $ \tau \in \mathbb{R} $ such that
$$
J'_{\lambda, k}(u_{0}) = \tau E'_{\lambda, k} (u_{0})
\quad \mbox{in} \quad  \left( W^{1,p(x)}(\mathbb{R}^{N}) \right)^{*},
$$
and so,
$$
0 = J'_{\lambda, k}(u_{0})  u_{0} = \tau E'_{\lambda, k} (u_{0})  u_{0}.
$$
If $ u_{0} \not\in \mathcal{M}_{\lambda, k}^{0} $, we must have $ E'_{\lambda, k} (u_{0})  u_{0} \neq 0 $. Hence, $ \tau = 0 $ and $ J'_{\lambda, k}(u_{0}) = 0 $ in  $ \left( W^{1,p(x)}(\mathbb{R}^{N}) \right)^{*} $, showing the lemma.
\end{pf}

\begin{lem}\label{Nehariempty}
Under the assumptions (\ref{p2}), \ref{H1} and \ref{H3}, we have that $ \mathcal{M}_{\lambda, k}^{0} = \emptyset $ for all $ k \in \mathbb{N}$ and $ 0 < \lambda < \Lambda_{1} = \Lambda_{1}(k) $, where
\begin{align}
\Lambda_{1} =  \frac{K^{-q_{+}}}{2  \Vert g_{k} \Vert_{\Theta(x)}} \left( \frac{r_{-} - p_{+}}{r_{+} - q_{-}}\right) \left[ \left(\frac{p_{-} - q_{+}}{r_{+} - p_{+}}\right) K^{- r_{+}} \right]^{\frac{p_{+} - q_{-}}{r_{-} - p_{+}}} .
\end{align}

\end{lem}
\begin{pf}
Arguing by contradiction, if the lemma does not hold, we have $ \mathcal{M}_{\lambda, k}^{0} \neq \emptyset $ for some $ \lambda_{0} \in (0, \Lambda_{1})$ and $k \in \mathbb{N}$.  Thereby, for $ u \in \mathcal{M}_{\lambda_{0}, k}^{0} $,
\begin{align*}
0 & =  E'_{\lambda_{0}, k} (u)u \\
  & =  \int p(x) (|\nabla u|^{p(x)} + |u|^{p(x)}) - \lambda_{0} \int q(x) g_{k}(x) |u|^{q(x)} - \int r(x) f_{k}(x) |u|^{r(x)} \\
    & \leq  p_{+}\int (|\nabla u|^{p(x)} + |u|^{p(x)}) - \lambda_{0} q_{-} \int  g_{k}(x) |u|^{q(x)} - r_{-}\int  f_{k}(x) |u|^{r(x)} \\
      & = \lambda_{0} (r_{-} - q_{-})\int  g_{k}(x)|u|^{q(x)} - (r_{1} - p_{+}) \int (|\nabla u|^{p(x)} + |u|^{p(x)}).
\end{align*}
By Propositions \ref{modular2}  and \ref{modular:interp}, H\"{o}lder's inequality and Sobolev embedding,
\begin{align}
\min\{ \Vert u \Vert^{p_{-}}, \Vert u \Vert^{p_{+}} \} \leq 2 \lambda_{0} \left( \frac{r_{-} - q_{-}}{r_{-} - p_{+}} \right) \Vert g_{k} \Vert_{\Theta(x)} K^{q_{+}}\max\{ \Vert u \Vert^{q_{-}}, \Vert u \Vert^{q_{+}} \}.\label{eql1}
\end{align}
Similarly,
\begin{align*}
0 & =  E'_{\lambda_{0}, k}(u) u \\
  & \geq  p_{-} \int (|\nabla u|^{p(x)} + |u|^{p(x)})   -  \lambda_{0} q_{+} \int g_{k}(x) |u|^{q(x)} - r_{+} \int f_{k}(x) |u|^{r(x)}\\
    &  = (p_{-} - q_{+}) \int (|\nabla u|^{p(x)} + |u|^{p(x)})  - (r_{+} - q_{+}) \int f_{k}(x) |u|^{r(x)}.
\end{align*}
Hence,
\begin{align}
\Big(\frac{p_{-} - q_{+}}{r_{+} - q_{+}}\Big) \min\{ \Vert u \Vert^{p_{-}}, \Vert u \Vert^{p_{+}} \} \leq K^{r_{+}}\max\{ \Vert u \Vert^{r_{-}}, \Vert u \Vert^{r_{+}} \}.\label{eql2}
\end{align}
If $ \Vert u \Vert \geq 1 $, it follows from (\ref{eql1}) that
\begin{align}
\Vert u \Vert \leq \left[ 2\lambda_{0} \left( \frac{r_{-} - q_{-}}{r_{-} - p_{+}}\right)  \Vert g_{k} \Vert_{\Theta(x)} K^{q_{+}} \right]^{\frac{1}{p_{-} - q_{+}}}.\label{eql3}
\end{align}
On the other hand, by (\ref{eql2}),
\begin{align}
\Vert u \Vert \geq \left[ \left(\frac{p_{-} - q_{+}}{r_{+} - q_{+}} \right) K^{-r_{+}}\right]^{\frac{1}{r_{+} - p_{-}}}. \label{eql4}
\end{align}
Combining (\ref{eql3}) and (\ref{eql4}), we derive that
\begin{align}
\lambda_{0} \geq \frac{K^{-q_{+}}}{2  \Vert g_{k} \Vert_{\Theta(x)}} \left( \frac{r_{-} - p_{+}}{r_{+} - q_{-}} \right) \left[ \left( \frac{p_{-} - q_{+}}{r_{+} - p_{+}} \right) K^{- r_{+}} \right]^{\frac{p_{-} - q_{+}}{r_{+} - p_{-}}}. \label{eql5}
\end{align}
Since
$$
0<\left(  \frac{p_{-} - q_{+}}{r_{+} - p_{+}} \right)K^{- r_{+}}  <1 \quad \mbox{and} \quad  \frac{p_{+} - q_{-}}{r_{-} - p_{+}} > \frac{p_{-} - q_{+}}{r_{+} - p_{-}},
$$
we deduce
$$
\lambda_0 \geq \frac{K^{-q_{+}}}{2  \Vert g_{k} \Vert_{\Theta(x)}} \left( \frac{r_{-} - p_{+}}{r_{+} - q_{-}} \right) \left[ \left(  \frac{p_{-} - q_{+}}{r_{+} - p_{+}} \right)K^{- r_{+}} \right]^{\frac{p_{+} - q_{-}}{r_{-} - p_{+}}},
$$
which is a contradiction.

Now if $ \Vert u \Vert < 1 $, we get from (\ref{eql1}),
\begin{align}
\Vert u \Vert \leq \left[ 2\lambda_{0} \left( \frac{r_{-} - q_{-}}{r_{-} - p_{+}} \right)  \Vert g_{k} \Vert_{\Theta(x)} K^{q_{+}} \right]^{\frac{1}{p_{+} - q_{-}}}.\label{eql6}
\end{align}
But by (\ref{eql2}),
\begin{align}
\Vert u \Vert \geq \left[ \left( \frac{p_{-} - q_{+}}{r_{+} - q_{+}} \right) K^{-r_{+}}\right]^{\frac{1}{r_{-} - p_{+}}}. \label{eql7}
\end{align}
Combining (\ref{eql6}) and (\ref{eql7}),
\begin{align}
\lambda_{0} \geq \frac{K^{-q_{+}}}{2  \Vert g_{k} \Vert_{\Theta(x)}} \left( \frac{r_{-} - p_{+}}{r_{+} - q_{-}} \right) \left[ \left(  \frac{p_{-} - q_{+}}{r_{+} - p_{+}} \right)K^{- r_{+}} \right]^{\frac{p_{+} - q_{-}}{r_{-} - p_{+}}}
\end{align}
so  a  new contradiction, finishing the proof.
\end{pf}

By Lemma~\ref{Nehariempty}, for $ 0 < \lambda < \Lambda_{1}$, we can write
$$
\mathcal{M}_{\lambda, k} = \mathcal{M}_{\lambda, k}^{+} \cup \mathcal{M}_{\lambda, k}^{-}.
$$
Therefore,  hereafter we will consider the following numbers

\begin{align*}
\alpha_{\lambda, k} = \inf_{u \in \mathcal{M}_{\lambda, k}} J_{\lambda, k}(u), \quad \alpha^{+}_{\lambda, k} = \inf_{u \in \mathcal{M}_{\lambda, k}^{+}} J_{\lambda, k}(u) \quad \mbox{ and } \quad \alpha^{-}_{\lambda, k} = \inf_{u \in \mathcal{M}_{\lambda, k}^{-}} J_{\lambda, k}(u).
\end{align*}

The next five lemmas establish important properties about the sets $\mathcal{M}_{\lambda, k}^{+}$ and $ \mathcal{M}_{\lambda, k}^{-}$.

\begin{lem}\label{Jlk-Mlk-plus}
Assume (\ref{p2}), \ref{H1}, \ref{H2} and \ref{H4}. If $ 0 < \lambda < \Lambda_{1}$,
then $ J_{\lambda, k}(u) < 0 $ for all $ u \in \mathcal{M}_{\lambda, k}^{+} $. \!\!Consequently, $ \alpha_{\lambda, k} \leq \alpha^{+}_{\lambda, k} < 0 .$
\end{lem}
\begin{pf}
Let $ u \in \mathcal{M}_{\lambda, k}^{+} $. Then, by definition of $ E'_{\lambda, k}(u)u $,
\begin{align*}
0 & < E'_{\lambda, k}(u) u \leq (r_{-} - q_{-})  \lambda  \int\!\! g_{k}(x) |u|^{q(x)} - (r_{-} - p_{+})\int \!(|\nabla u|^{p(x)} + |u|^{p(x)}), 
\end{align*}
from where it follows
\begin{equation} \label{eql8}
 \lambda \int g_{k}(x) |u|^{q(x)} > \Big(\frac{r_{-} - p_{+}}{r_{-} - q_{-}}\Big) \int (|\nabla u|^{p(x)} + |u|^{p(x)}).
\end{equation}
By definition of $ J_{\lambda, k}(u) $,
\begin{align*}
J_{\lambda, k}(u) & \leq \frac{1}{p_{-}} \int (|\nabla u|^{p(x)} + |u|^{p(x)}) - \frac{\lambda}{q_{+}}\int g_{k}(x) |u|^{q(x)} - \frac{1}{r_{+}}\int f_{k}(x) |u|^{r(x)} \\
 & = \left( \frac{1}{p_{-}} - \frac{1}{r_{+}} \right) \int (|\nabla u|^{p(x)} + |u|^{p(x)}) - \lambda \left( \frac{1}{q_{+}} - \frac{1}{r_{+}} \right) \int g_{k}(x) |u|^{q(x)}.
\end{align*}
By (\ref{eql8}) and \ref{H4},
\begin{align*}
J_{\lambda, k}(u) & \leq \left[ \frac{1}{p_{-}} - \frac{1}{r_{+}} - \left( \frac{1}{q_{+}} - \frac{1}{r_{+}} \right) \left( \frac{r_{-} - p_{+}}{r_{-} - q_{-}}\right) \right] \int (|\nabla u|^{p(x)} + |u|^{p(x)}) \\
  & = \left[ \frac{r_{+} - p_{-}}{p_{-}r_{+}} - \Big(\frac{r_{+} - q_{+}}{q_{+}r_{+}}\Big) \cdot  \Big(\frac{r_{-} - p_{+}}{r_{-} - q_{-}}\Big) \right] \int (|\nabla u|^{p(x)} + |u|^{p(x)}) \\
    & = \frac{(r_{+} - p_{-})}{r_+} \left[ \frac{1}{p_{-}} - \frac{1}{q_{+}} \cdot \Big( \frac{r_{+} - q_{+}}{r_{+} - p_{-}}\Big) \cdot \Big( \frac{r_{-} - p_{+}}{r_{-} - q_{-}}\Big) \right] \int (|\nabla u|^{p(x)} + |u|^{p(x)}) \\
     & < 0.
\end{align*}
\end{pf}

\begin{lem}\label{desMlk}
We have the following inequalities
\begin{enumerate}
\item[(i)] $ \int g_{k}(x) |u|^{q(x)} > 0 $ for each $ u \in \mathcal{M}_{\lambda, k}^{+}$;
\item[(ii)] $ \Vert u \Vert < \left[ 2 \Big(\frac{r_{-} - q_{-}}{r_{-} - p_{+}}\Big) K^{q_{+}}  \right]^{1/(p_{-} - q_{+})}\max\left\{ \left(\lambda \Vert g_{k} \Vert_{\Theta(x)} \right)^{\frac{1}{p_{+} - q_{-}}}, \left(\lambda \Vert g_{k} \Vert_{\Theta(x)}  \right)^{\frac{1}{p_{-} - q_{+}}}\right\} $ for each $ u \in \mathcal{M}_{\lambda, k}^{+}$;
\item[(iii)] $ \Vert u \Vert > \left[ \Big(\frac{p_{-} - q_{+}}{r_{+} - q_{+}} \Big) K^{-r_{+}}\right]^{\frac{1}{r_{+} - p_{-}}} $ for each $ u \in \mathcal{M}_{\lambda, k}^{-}$.
\end{enumerate}
\end{lem}

\begin{pf} \\
\noindent (i) An immediate consequence of (\ref{eql8}).  \\
\noindent (ii) Similarly  to the proof of Lemma~\ref{Nehariempty},
\begin{align*}
\min\{ \Vert u \Vert^{p_{-}}, \Vert u \Vert^{p_{+}} \} < 2 \lambda \Big(\frac{r_{-} - q_{-}}{r_{-} - p_{+}}\Big) \Vert g_{k} \Vert_{\Theta(x)} K^{q_{+}}\max\{ \Vert u \Vert^{q_{-}}, \Vert u \Vert^{q_{+}} \}.
\end{align*}
If $ \Vert u \Vert < 1 $, the above inequality gives
\begin{align*}
\Vert u \Vert \leq \left[ 2 \lambda \Big(\frac{r_{-} - q_{-}}{r_{-} - p_{+}}\Big) \Vert g_{k} \Vert_{\Theta(x)} K^{q_{+}} \right]^{\frac{1}{p_{+} - q_{-}}}.
\end{align*}
Now, if $ \Vert u \Vert \geq 1 $, we will get
\begin{align*}
\Vert u \Vert \leq \left[ 2 \lambda \Big(\frac{r_{-} - q_{-}}{r_{-} - p_{+}}\Big) \Vert g_{k} \Vert_{\Theta(x)} K^{q_{+}} \right]^{\frac{1}{p_{-} - q_{+}}},
\end{align*}
showing (ii).

\bigskip
\noindent (iii) Let $ u \in \mathcal{M}_{\lambda, k}^{-} $. Similarly to the proof of Lemma~\ref{Nehariempty},
$$
\Big(\frac{p_{-} - q_{+}}{r_{+} - q_{+}}\Big) \min\{ \Vert u \Vert^{p_{-}}, \Vert u \Vert^{p_{+}} \} < K^{r_{+}}\max\{ \Vert u \Vert^{r_{-}}, \Vert u \Vert^{r_{+}} \}.
$$
If $ \Vert u \Vert < 1 $,
\begin{equation} \label{DES1}
\Vert u \Vert > \left[ \Big(\frac{p_{-} - q_{+}}{r_{+} - q_{+}}\Big)K^{-r_{+}} \right]^{\frac{1}{r_{-} - p_{+}}},
\end{equation}
and for $ \Vert u \Vert \geq 1 $,
\begin{equation} \label{DES2}
\Vert u \Vert > \left[ \Big(\frac{p_{-} - q_{+}}{r_{+} - q_{+}}\Big)K^{-r_{+}} \right]^{\frac{1}{r_{+} - p_{-}}}.
\end{equation}
Thus, the last two inequalities imply that (iii) hold.
\end{pf}

\begin{lem}\label{Jlk-Mlk-menos}
Assume that $ 0 < \lambda < \frac{q_{-}}{p_{+}}\Lambda_{1}$  and ~\ref{H1}. Then there exists a positive constant $ d_{1} = d_{1}(p_{\pm}, q_{\pm}, r_{\pm},  K, \Vert g_{k} \Vert_{\Theta(x)})$ such that $ J_{\lambda, k}(u) > 0 $ for each $ u \in \mathcal{M}_{\lambda, k}^{-} $.
\end{lem}

\begin{pf}
Let $ u \in \mathcal{M}_{\lambda, k}^{-} $. Then, using the definitions of $ J_{\lambda, k} $  and $ \mathcal{M}_{\lambda, k} $, we can write
\begin{align*}
J_{\lambda, k}(u) & \geq \left( \frac{1}{p_{+}} - \frac{1}{r_{-}} \right) \int (|\nabla u|^{p(x)} + |u|^{p(x)}) - \lambda \left( \frac{1}{q_{-}} - \frac{1}{r_{-}} \right) \int g_{k}(x) |u|^{q(x)} \\
& \geq \left( \frac{1}{p_{+}} - \frac{1}{r_{-}} \right)  \min\{ \Vert u \Vert^{p_{-}}, \Vert u \Vert^{p_{+}} \} \\
 & \qquad- 2 \lambda \left( \frac{1}{q_{-}} - \frac{1}{r_{-}} \right) \Vert g_{k} \Vert_{\Theta(x)} K^{q_{+}} \max\{ \Vert u \Vert^{q_{-}}, \Vert u \Vert^{q_{+}} \}.
\end{align*}
If $ \Vert u \Vert < 1 $, it follows that
\begin{align*}
J_{\lambda, k}(u) & \geq \left( \frac{1}{p_{+}} - \frac{1}{r_{-}} \right) \Vert u \Vert^{p_{+}} - 2 \lambda \left( \frac{1}{q_{-}} - \frac{1}{r_{-}} \right) \Vert g_{k} \Vert_{\Theta(x)} K^{q_{+}} \Vert u \Vert^{q_{-}}\\
  & = \Vert u \Vert^{q_{-}} \left[ \left( \frac{1}{p_{+}} - \frac{1}{r_{-}} \right) \Vert u \Vert^{p_{+ } - q_{-}} - 2 \lambda \left( \frac{1}{q_{-}} - \frac{1}{r_{-}} \right) \Vert g_{k} \Vert_{\Theta(x)} K^{q_{+}} \right].
\end{align*}
Thereby, by Lemma~\ref{desMlk}~(iii),
\begin{align*}
J_{\lambda, k} (u) & > \left[ \left(\frac{p_{-} - q_{+}}{r_{+} - q_{+}} \right) K^{-r_{+}}  \right]^{\frac{q_{-} }{r_{+} - p_{-}}} \Biggl\{ \left( \frac{1}{p_{+}} - \frac{1}{r_{-}} \right)\left[ \left(\frac{p_{-} - q_{+}}{r_{+} - q_{+}}\right)K^{-r_{+}} \right]^{\frac{p_{-} - q_{+}}{r_{+} - p_{-}}}  \\
&  \qquad\quad - 2 \lambda \left( \frac{1}{q_{-}} - \frac{1}{r_{-}} \right) \Vert g_{k} \Vert_{\Theta(x)} K^{q_{+}} \Biggr\}=d_1.
\end{align*}
Similarly, if $ \Vert u \Vert \geq 1 $,
\begin{align*}
J_{\lambda, k} (u) & > \left[ \left( \frac{p_{-} - q_{+}}{r_{+} - q_{+}} \right) K^{-r_{+}}  \right]^{\frac{q_{+}}{r_{+} - p_{-}}} \Biggl\{ \left( \frac{1}{p_{+}} - \frac{1}{r_{-}} \right)\left[ \left(\frac{p_{-} - q_{+}}{r_{+} - q_{+}}\right)K^{-r_{+}} \right]^{\frac{p_{-} - q_{+}}{r_{+} - p_{-}}}   \\
&  \qquad\quad  - 2 \lambda \left( \frac{1}{q_{-}} - \frac{1}{r_{-}} \right) \Vert g_{k} \Vert_{\Theta(x)} K^{q_{+}} \Biggr\}=d_1.
\end{align*}
From the above estimates, the lemma follows if $ 0 < \lambda <\frac{q_{-}}{p_{+}} \Lambda_{1} $.

\end{pf}


The lemma below  is crucial in our arguments because it shows a condition for the existence of exactly two nontrivial zeroes for a special class of functions.

\begin{lem}\label{lemAbst}
Let $ g_{i} : [0, + \infty) \to [0, +\infty) $, $ i \in \{1,2,3\} $,  be  increasing continuous functions, with $ g_{i}(0) =0 $   verifying the following conditions:
\begin{enumerate}
\item[(i)]  $\displaystyle \lim_{t \to 0^+ }\frac{g_{3}(t)}{g_{1}(t)} = 0 $;
\item[(ii)] $\displaystyle \lim_{t \to + \infty} g_{2}(t) = +\infty $;

\item[(iii)] $\displaystyle \lim_{t \to 0^+} \frac{g_{1}(t) - g_{3}(t)}{g_{2}(t)} = 0$;

\item[(iv)] the function $ \phi = g_{1} - g_{3} $ has only one maximum point and $ \phi(t) \to -\infty $ as $ t \to +\infty $.

\end{enumerate}
Suppose  there exists $ \tilde{t} \in (0, t_{\max}) $ with $ \phi(t_{\max}) = \max_{t \geq 0} \phi(t) $ such that $ \frac{g_{1} - g_{3}}{g_{2}} $ is increasing on $(0, \tilde{t} ) $.  Then, there is $ \lambda_{*} > 0 $ such that  $ \psi = g_{1} - \lambda g_{2} - g_{3} $  has only two nontrivial zeros for all $ 0 <\lambda < \lambda_{*} $.
\end{lem}
\begin{pf}
From  ($i$), it is clear that $ \phi(t) > 0 $ for all $ t > 0 $ sufficiently small. Since $ \frac{g_{1} - g_{3}}{g_{2}} $ is positive and increasing in the interval $(0, \tilde{t} ) $, for each $ 0  < \lambda < \phi(\tilde{t}) $ there is unique $ t_{\lambda} \in (0, \tilde{t}) $ such that
\begin{align*}
\lambda = \frac{g_{1}(t_{\lambda})- g_{3}(t_{\lambda})}{g_{2}(t_{\lambda})}.
\end{align*}
Then, by hypothesis that $ \frac{g_{1} - g_{3}}{g_{2}} $ is increasing on $(0, \tilde{t} ) $, we derive
\begin{align*}
\lambda g_{2}(t) < g_{1}(t)- g_{3}(t)\quad \text{for all } t \in (t_{\lambda}, \tilde{t}).
\end{align*}
Now, fix  $ \lambda^* > 0 $ such that
\begin{align*}
\lambda g_{2}(t) < g_{1}(t)- g_{3}(t)\quad \text{for all } t \in (t_{\lambda}, t_{\max}) \quad \mbox{and} \quad \lambda \in (0,\lambda^{*}).
\end{align*}


Since $ \phi $ is decreasing in the interval $ (t_{\max}, \infty)$, $ g_{2} $ is increasing and  $ g_{2}(t) \to \infty $ as $ t \to \infty $, there is a unique number $ t_{1} > t_{\max} $ such that
\begin{align*}
\lambda g_{2}(t_{1}) = \phi(t_{1}).
\end{align*}
Therefore, $t_\lambda$ and $t_1$ are the unique nontrivial zeros of $\psi$ for $\lambda \in (0,\lambda^*)$.
\end{pf}

With the help of Lemma~\ref{lemAbst}, we get the following result, which is similar  to the constant case, see \cite{BrownWu} and \cite{Fan2008}.
\begin{lem}\label{Max-Nehari}
For each $ u \in W^{1, p(x)}(\mathbb{R}^{N})\setminus \{0\} $, we have the following:
\begin{enumerate}
\item[(i)] if $ \int g_{k}(x) |u|^{q(x)} = 0 $, then there exists a unique positive number $ t^{-} = t^{-}(u)
 $ such that
 \[
 t^{-}u \in \mathcal{M}^{-}_{\lambda, k} \quad \mbox{ and } \quad J_{\lambda, k}(t^{-}u) = \sup_{t \geq 0}J_{\lambda, k}(t u);
 \]
\item[(ii)] if $ 0 < \lambda < \Lambda_{1}$ and $ \int g_{k}(x) |u|^{q(x)} > 0 $, then there exist $ t^{*} > 0 $ and unique positive numbers $ t^{+} = t^{+}(u) < t^{-} = t^{-}(u) $ such that $ t^{+} u \in \mathcal{M}^{+}_{\lambda, k} $, $ t^{-} u \in \mathcal{M}^{-}_{\lambda, k} $ and
    \[
    J_{\lambda, k}(t^{+} u) = \inf_{0 \leq t \leq t^{*}} J_{\lambda, k}(t u), \quad J_{\lambda, k}(t^{-} u) = \sup_{t \geq t^{*}} J_{\lambda, k}(t u).
    \]
\end{enumerate}
\end{lem}
\begin{pf}
By direct calculations, we see that
\begin{align*}
E'_{\lambda, k}(tu) tu =  t \frac{d}{dt}(J_{\lambda, k} (tu)) + t^{2} \frac{d^{2}}{dt^{2}}(J_{\lambda, k} (tu)).
\end{align*}
Thus, if $ t = \bar t $ is a critical point of $ J_{\lambda, k}(tu) $,
\begin{align}\label{eq:der_Jtu}
E'_{\lambda, k}(\bar t u) \bar t u = \bar {t\,}^{2} \left. \frac{d^{2}}{dt^{2}}(J_{\lambda, k} (tu))\right|_{t = \bar t}.
\end{align}
Using \eqref{eq:der_Jtu} and the same ideas of the proof of Lemma~3.6 of \cite{Fan2008}, we get the item~(i).

To prove item (ii), we will apply the Lemma~\ref{lemAbst} with the functions:
\begin{align*}
g_{1}(t) & = \int t^{p(x) - 1}(|\nabla u|^{p(x)} + |u|^{p(x)});\\
g_{2}(t) & = \int t^{q(x) - 1}g_{k}(x) |u|^{q(x)};
\end{align*}
and
\begin{align*}
g_{3}(t) = \int t^{r(x) - 1} f_{k}(x)|u|^{r(x)}.
\end{align*}
The reader is invited to check that $ g_{1}, g_{2} $ and $ g_{3} $ satisfy the conditions of Lemma~\ref{lemAbst}, and so, the function $ \psi(t) = g_{1}(t) - \lambda g_{2}(t) - g_{3}(t) = J'_{\lambda, k}(tu) u $ has only two nontrivial zeros, $ t^{+} < t^{-} $. Let $ \varphi(t) = J_{\lambda, k}(tu) $ on $ [0, \infty) $. Then, it is clear that $ \varphi(0) = 0 $ and $ \varphi(t) $ is negative if $ t > 0 $ is small, implying that $\varphi$  has a local minimum in $ t = t^{+} $.  Consequently,
\begin{align*}
E'_{\lambda, k}(t^{+}u) t^{+}u > 0,
\end{align*}
from where it follows that $ t^{+}u \in \mathcal{M}^{+}_{\lambda, k} $. Since $ t^{+} $ and $ t^{-} $ are the unique critical points of $ \varphi $, we deduce that  $ \varphi $ has a global maximum in $ t = t^{-} $, thus
\begin{align*}
E'_{\lambda, k}(t^{+}u) t^{+}u < 0.
\end{align*}
and $ t^{-}u \in \mathcal{M}^{-}_{\lambda, k} $.  Using the Lemmas~\ref{Jlk-Mlk-plus} and \ref{Jlk-Mlk-menos}, it follows that $ J_{\lambda, k}(t^{+}u) < 0 $ and $ J_{\lambda, k}(t^{-}u) > 0 $. Let $ t_{*} > 0 $ be the unique zero of  $\varphi $ in $ (t^{+}, t^{-})$. Then is clear that
\begin{align*}
J_{\lambda, k}(t^{+}u) = \inf_{0 \leq t \leq t_{*}} J_{\lambda, k}(tu) \qquad \mbox{ and } \qquad J_{\lambda, k}(t^{-}u) = \max_{t \geq t_{*}} J_{\lambda, k}(tu).
\end{align*}
\end{pf}

\begin{lem}\label{ltdalem}
Assume that $ g $ satisfies \ref{H1} and let $ \{ u_{n} \} $ be a $ (PS)_{d} $ sequence in $ W^{1, p(x)}(\mathbb{R}^{N}) $ for $ J_{\lambda, k} $.  Then $ \{u_{n}\} $ is bounded in $ W^{1,p(x)}(\mathbb{R}^{N}) $.
\end{lem}
\begin{pf}
It is clear that
\begin{align*}
J_{\lambda, k}(u_{n}) - \frac{1}{r_{-}} J'_{\lambda, k}(u_{n}) u_{n}  
   \geq& \left(  \frac{1}{p_{+}} - \frac{1}{r_{-}} \right) \int \left( |\nabla u_{n}|^{p(x)} +  |u_{n}|^{p(x)}\right) \\
& \quad + \lambda \left( \frac{1}{r_{-}} - \frac{1}{q_{-}} \right) \int g_{k}(x) |u_{n}|^{q(x)}.
\end{align*}

Assume that $ \Vert u_ {n} \Vert \geq  1 $ for some $ n \in \mathbb{N} $. Then, by H\"{o}lder's inequality and Sobolev embedding, we derive the inequality
\begin{align*}
d + 1 + \Vert u_{n} \Vert \geq \left(  \frac{1}{p_{+}} - \frac{1}{r_{-}} \right) \Vert u_{n} \Vert^{p_{-}} - \lambda \left( \frac{1}{q_{-}} - \frac{1}{r_{-}} \right) \Vert g_{k} \Vert_{\Theta(x)} K^{q_{+}} \Vert u_{n} \Vert^{q_{+}}.
\end{align*}
Since $ 1 < q_{+} < p_ {-}$, the last inequality yields $ \{u_{n} \} $ is bounded in $ W^{1,p(x)}(\mathbb{R}^{N}) $.
\end{pf}

Now, combining standard arguments with the boundedness of $\{u_n\}$ and Sobolev imbedding (see \cite{AlvesBarreiro2}), we have the below result.
\begin{thm}\label{conv}
Assume that $ g $ satisfies \ref{H1}. If $ \{ u_{n} \} $ is a sequence in $  W^{1,p(x)}(\mathbb{R}^{N}) $ such that $ u_{n} \rightharpoonup u $ in $  W^{1,p(x)}(\mathbb{R}^{N}) $ and $ J'_{\lambda, k}(u_{n}) \to 0 $ as $ n \to \infty $, then for some subsequence, still denoted by $\{u_{n}\}$, $ \nabla u_{n}(x) \to \nabla u(x) $ a.e. in $ \mathbb{R}^{N} $  and $ J'_{\lambda, k}(u) = 0 $.
\end{thm}

The next theorem is a compactness result  on Nehari manifolds. The case for constant exponent is due to Alves \cite{alves05}. 

\begin{thm}\label{TeoComp}
Suppose that (\ref{p2}) holds and let $ \{u_{n}\} \subset \mathcal{M}_{\infty} $ be a sequence with $ J_{\infty}(u_{n}) \to c_{\infty} $. Then,
\begin{description}
\item[I.] $ u_{n} \to u $ in $ W^{1,p(x)}(\mathbb{R}^{N}) $,

or

\item[II.] There is $ \{y_{n}\} \subset \mathbb{Z}^{N}$ with $|y_n| \to +\infty$ and $w \in W^{1,p(x)}(\mathbb{R}^{N}) $ such that $ w_{n} (x) = u_{n}(x + y_{n}) \to w $ in $ W^{1,p(x)}(\mathbb{R}^{N}) $ and $J_{\infty}(w) = c_{\infty}$.
\end{description}
\end{thm}
\noindent \textbf{Proof.} Similarly to Corollary~\ref{ltdalem}, there is $ u \in W^{1,p(x)}(\mathbb{R}^{N}) $ and a subsequence of $ \{ u_{n} \}$, still denoted by itself, such that $u_n  \rightharpoonup u $  in $ W^{1,p(x)}(\mathbb{R}^{N})$. Applying  Ekeland's variational principle, we  can assume that 
\begin{align} \label{eq1}
J'_{\infty}(u_{n}) - \tau_{n} E'_{\infty}(u_{n}) = o_{n}(1),
\end{align}
where $ (\tau_{n}) \subset \mathbb{R} $ and $ E_{\infty}(w) = J'_{\infty}(w)  w $, for any $ w \in W^{1,p(x)}(\mathbb{R}^{N}) $.

Since $ \{u_{n}\} \subset \mathcal{M}_{\infty} $, (\ref{eq1}) leads to
\[
\tau_{n} E'_{\infty}(u_{n})  u_{n} = o_{n}(1).
\]
Next, we will show that there exists $\eta  > 0$ such that
\begin{align} \label{ldtaMoo}
|E'_{\infty}(u_{n})u_{n}| > \eta \,\,\, \forall n \in \mathbb{N}.
\end{align}
Indeed, first we claim that there exists $\eta_{0}  > 0 $ satisfying
\[
\Vert u \Vert > \eta_{0} \quad \mbox{ for any } u \in \mathcal{M}_{\infty}.
\]
Suppose by contradiction that the claim is false. Then, there is $ \{ v_n \} \subset \mathcal{M}_{\infty} $ such that $ \Vert v_n \Vert \to 0 $ as $ n \to \infty $. Since $ \{ v_n \} \subset \mathcal{M}_{\infty} $, we derive
\[
\int \left( |\nabla v_{n}|^{p(x)} + |v_{n}|^{p(x)} \right) = \int |v_{n}|^{r(x)}.
\]
On the other hand, using the fact that $ \Vert v_{n} \Vert < 1 $ for $ n $ large enough, it follows from Propositions~\ref{modular} and \ref{modular2},
\begin{align*}
\Vert v_{n} \Vert^{p_{+}} \leq C\max\{ \Vert v_{n} \Vert^{r_{-}}, \Vert v_{n} \Vert^{r_{+}} \}=C\|v_n\|^{r_-},
\end{align*}
leading to
$$
\left(\frac{1}{C}\right)^{\frac{1}{r_- - p_+}}\leq \|v_n\|,
$$
which is absurd. Therefore, by Proposition~\ref{modular2}, there is $ \varsigma > 0 $ such that
\[ \rho_{1}(u) \geq  \varsigma \quad u \in \mathcal{M}_{\infty}.
\]
By definition of $ E_{\infty}(u) $,
\begin{align*}
E'_{\infty}(u_n)u_n & \leq p_{+} \int \left( |\nabla u_{n}|^{p(x)} + |u_{n}|^{p(x)} \right) - r_{-} \int |u_{n}|^{r(x)} \\
& = (p_{+} - r_{-}) \int \left( |\nabla u_{n}|^{p(x)} + |u_{n}|^{p(x)} \right) =(p_{+} - r_{-}) \rho_1(u_n)< (p_{+} - r_{-}) \varsigma,
\end{align*}
proving (\ref{ldtaMoo}).  Now, combining (\ref{eq1}) and (\ref{ldtaMoo}), we see that $ \tau_{n} \to 0 $, and so,
\[
J_{\infty}(u_{n}) \to c_{\infty} \,\,\, \mbox{and} \,\,\, J'_{\infty}(u_{n}) \to 0.
\]

Next, we will study the following  possibilities: $ u \neq 0 $ or $ u = 0 $.

\vspace{0.5 cm}

\noindent \textbf{Case 1:} $ u \neq 0 $.

\vspace{0.5 cm}

Similarly to Theorem~\ref{conv},  it follows that the below limits are valid for some subsequence:
\begin{itemize}
    \item $ u_{n}(x) \to u(x)$ \,\, and \,\, $\nabla u_{n}(x) \to \nabla u(x) $ a.e. in $ \mathbb{R}^{N}, $
    \item $ \displaystyle \int |\nabla u_{n}(x)|^{p(x)-2}\nabla u_{n}(x)\nabla v \to \int |\nabla u(x)|^{p(x)-2}\nabla u(x)\nabla v$,
    \item $  \displaystyle \int |u_{n}|^{p(x)-2} u_{n} v \to \int |u|^{p(x)-2} u_{n} v$,
\end{itemize}
and
\begin{itemize}
    \item $ \displaystyle \int |u_{n}|^{r(x) -2} u_{n} v \to \int |u|^{r(x) -2} u v  $
\end{itemize}
for any $ v \in W^{1,p(x)}(\mathbb{R}^{N})$. Consequently, $ u $ is critical point of  $ J_{\infty} $.
By Fatou's Lemma , it is easy to check that
\begin{align*}
c_{\infty} \leq & J_{\infty}(u) = J_{\infty}(u) - \frac{1}{r_{-}} J'_{\infty}(u) u \\
= &\int \left( \frac{1}{p(x)} - \frac{1}{r_{-}} \right) \left( | \nabla u |^{p(x)} + |u|^{p(x)} \right) + \int \left( \frac{1}{r_{-}} - \frac{1}{r(x)} \right)  |u|^{r(x)} \\
\leq & \liminf_{n \to \infty}\left\{ \int \left( \frac{1}{p(x)} - \frac{1}{r_{-}} \right) \left( | \nabla u_{n} |^{p(x)} + |u_{n}|^{p(x)} \right) \right.\\
&\qquad \left. + \int \left( \frac{1}{r_{-}} - \frac{1}{r(x)} \right)  |u_{n}|^{r(x)}\right\} \\
= & \liminf_{n \to \infty} \left\{ J_{\infty}(u_{n}) - \frac{1}{r_{-}} J'_{\infty}(u_{n}) u_{n} \right\} = \, c_{\infty} .\\
\end{align*}
Hence,
\begin{align*}
\lim_{n \to \infty} \int \left( | \nabla u_{n} |^{p(x)} + |u_{n}|^{p(x)} \right) = \int \left( | \nabla u |^{p(x)} + |u|^{p(x)} \right),
\end{align*}
implying that $ u_{n} \to u $ in $ W^{1,p(x)}(\mathbb{R}^{N})$.

\vspace{0.5 cm}

\noindent \textbf{Case 2:} $ u = 0 $.

\vspace{0.5 cm}

In this case, we claim that there are $R, \xi>0$ and $ \{ y_{n} \} \subset \mathbb{R}^{N} $ satisfying
\begin{align}\label{lionsfalse}
\limsup_{n \to \infty} \int_{B_{R}(y_{n})} |u_{n}|^{p(x)} \geq \xi.
\end{align}
If the claim is false, we must have
\begin{align*}
\limsup_{n \to \infty} \sup_{y \in \mathbb{R}^{N}} \int_{B_{R}(y)} |u_{n}|^{p(x)} = 0.
\end{align*}
Thus,  by a  Lions  type result  for variable exponent proved in \cite[Lemma~3.1]{Fan2001},
$$
u_{n} \to 0  \mbox{ in }  L^{s(x)}(\mathbb{R}^{N}),
$$
for any $ s \in C(\mathbb{R}^{N}) $ with $ p \ll s \ll p^{*} $.

 Recalling $ J'_{\infty}(u_{n}) u_{n} = o_{n}(1) $, the last limit yields
\[
\int \left( | \nabla u_{n} |^{p(x)} + |u_{n}|^{p(x)} \right) = o_{n}(1),
\]
or equivalently
$$
u_n \to 0 \,\,\, \mbox{in} \,\,\, W^{1,p(x)}(\mathbb{R}^{N}),
$$
leading to $c_{\infty} = 0 $, which is absurd. This way, (\ref{lionsfalse}) is true. By a routine argument, we can assume that $ {y}_{n} \in \mathbb{Z}^{N} $ and $ |y_{n}| \to \infty $ as $ n \to \infty $.
Setting
$$
w_{n}(x) = u_{n}(x + {y}_{n}),
$$
and using the fact that $ p $ and $ r $ are $ \mathbb{Z}^{N}$-periodic, a change of variable gives
$$
J_{\infty}(w_{n}) = J_{\infty}(u_{n}) \,\,\, \mbox{and} \,\,\, \|J'_{\infty}(w_{n})\| = \|J'_{\infty}(u_{n})\|,
$$
showing that $ \{ w_{n} \} $ is a sequence $(PS)_{c_{\infty}} $ for $ J_{\infty} $.  If $ w \in  W^{1,p(x)}(\mathbb{R}^{N}) $ denotes the weak limit of $ \{ w_{n} \} $, from (\ref{lionsfalse}),
$$
\int_{B_{{R}}(0)} |w|^{p(x)} \geq \xi,
$$
showing that $ w \neq 0 $.

Repeating the same argument of the first case for the sequence $ \{w_{n}\} $, we deduce that $ w_{n} \to w $ in $ W^{1,p(x)}(\mathbb{R}^{N}) $, $ w \in \mathcal{M}_{\infty} $  and $ J_{\infty}(w) = c_{\infty} $.  \fim

\vspace{0.5 cm}
Our next result will be very useful in the study of the compactness of some functionals.

\begin{lem}\label{J-ltda}
Let $u \in W^{1,(x)}(\mathbb{R}^{N})$ be a nontrivial critical point of $J_{\lambda, k}$.  Then, there exists a constant $ M = M(k) > 0 $, which is independent of $ \lambda $, such that
\begin{align*}
J_{\lambda, k}(u) \geq - M\left( \lambda^{\frac{p_{+}}{p_{+} - q_{-}}} + \lambda^{\frac{p_{-}}{p_{-} - q_{+}}}  \right).
\end{align*}
\end{lem}
\begin{pf}
By hypothesis, $ J'_{\lambda, k}(u)u = 0 $. Arguing as in the proof of  Lemma~\ref{Jlk-Mlk-menos},
if $ \Vert u \Vert \geq 1 $, then
\begin{align*}
J_{\lambda, k}(u) \geq \left( \frac{1}{p_{+}} -  \frac{1}{r_{-}} \right) \Vert u \Vert^{p_{-}} - \left( \frac{1}{q_{-}} -  \frac{1}{r_{-}} \right) 2\lambda \Vert g_{k} \Vert_{\Theta(x)} K^{q_{+}} \Vert u \Vert^{q_{+}}.
\end{align*}
Applying Young's inequality with $ p_{1} = \frac{p_{-}}{q_{+}} $ and $ p_{2} = \frac{p_{-}}{p_{-} - q_{+}} $, we obtain
\begin{align*}
J_{\lambda, k}(u) &\geq  \left( \frac{1}{p_{+}} -  \frac{1}{r_{-}} \right) \Vert u \Vert^{p_{-}} -\epsilon \left( \frac{1}{q_{-}} -  \frac{1}{r_{-}} \right)   \Vert u \Vert^{p_{-}} \\
& \quad- \left( \frac{1}{q_{-}} -  \frac{1}{r_{-}} \right)C_{1}(\epsilon) \left( 2\lambda  \Vert g_{k} \Vert_{\Theta(x)}K^{q_{+}} \right)^{\frac{p_{-}}{p_{-} -  q_{+}}}
\end{align*}
where $ C_{1}(\epsilon) = \frac{p_{-} - q_{+}}{p_{-}} \left( \frac{q_{+}}{\epsilon p_{-}} \right)^{\frac{q_{+}}{p_{-} - q_{+}}} $. Choosing $ \epsilon = \left( \frac{1}{q_{-}} -  \frac{1}{r_{-}} \right)^{-1} \left( \frac{1}{p_{+}} -  \frac{1}{r_{-}} \right) $, we get
\begin{align*}
J_{\lambda, k}(u) \geq - \left( \frac{1}{q_{-}} -  \frac{1}{r_{-}} \right)C_{1} (\epsilon) \left( 2\lambda  \Vert g_{k} \Vert_{\Theta(x)}K^{q_{+}} \right)^{\frac{p_{-}}{p_{-} -  q_{+}}}.
\end{align*}
Analogously, if $ \Vert u \Vert < 1 $, we will get
\begin{align*}
J_{\lambda, k}(u) \geq - \left( \frac{1}{q_{-}} -  \frac{1}{r_{-}} \right)C_{2} (\epsilon) \left(2 \lambda  \Vert g_{k} \Vert_{\Theta(x)} K^{q_{+}} \right)^{\frac{p_{+}}{p_{+} -  q_{-}}},
\end{align*}
where $ C_{2}(\epsilon) = \frac{p_{+} - q_{-}}{p_{+}} \left( \frac{q_{-}}{\epsilon p_{+}} \right)^{\frac{q_{-}}{p_{+} - q_{-}}} $.

Therefore,
\begin{align*}
J_{\lambda, k}(u) \geq - M\left( \lambda^{\frac{p_{+}}{p_{+} - q_{-}}} + \lambda^{\frac{p_{-}}{p_{-} - q_{+}}}  \right)
\end{align*}
with
$$
M = \left( \frac{1}{q_{-}} -  \frac{1}{r_{-}} \right)(2K^{q_+})^{\frac{p_{-}}{p_{-} - p_{+}}} \max\left\{C_{1}(\epsilon)\Vert g_{k} \Vert_{\Theta(x)}^{\frac{p_{-}}{p_{-} -  q_{+}}}, C_{2}(\epsilon)\Vert g_{k} \Vert_{\Theta(x)}^{\frac{p_{+}}{p_{+} -  q_{-}}}  \right\}.
$$
\end{pf}

The next result is an important step to prove the existence of solutions, because it establishes the behavior of the $(PS)$ sequences of functional $ J_{\lambda, k} $.
\begin{lem}
Let $ \{v_{n} \} $ be a $(PS)_{d}$ sequence   for functional  $ J_{\lambda, k} $ with $ v_{n} \rightharpoonup v $ in $ W^{1,p(x)}(\mathbb{R}^{N}) $ . Then,
\begin{align}
J_{\lambda, k} (v_{n}) - J_{0, k} (w_{n}) - J_{\lambda, k}(v) = o_{n}(1) \label{J-J0-on1}
\end{align}
and
\begin{align}
\Vert J'_{\lambda, k} (v_{n}) - J'_{0, k} (w_{n}) - J'_{\lambda, k}(v) \Vert = o_{n}(1),\label{J'-J0'-on1}
\end{align}
where $ w_{n} = v_{n} - v $.
\end{lem}
\begin{pf}
Similarly to Theorem~\ref{conv}, the below limits occur
\[
\nabla v_{n}(x) \to \nabla v(x) \,\,\, \mbox{and} \,\,\, v_{n}(x) \to v(x) \,\,\, \mbox{ a.e. in} \,\,\, \mathbb{R}^{N}.
\]
Then, by Proposition \ref{Brezis-Lieb-1},
\begin{align*}
J_{\lambda, k}(v_{n})  = J_{0,k}(w_{n}) + J_{\lambda, k}(v) + o_{n}(1), \label{eq7}
\end{align*}
showing~(\ref{J-J0-on1}). The equality (\ref{J'-J0'-on1}) follows from Propositions \ref{Brezis-Lieb-2} and \ref{Brezis-Lieb-3}.
\end{pf}

The proof of the next result follows the same steps found in \cite{Lin12} and \cite{Miotto}, and so, it will be omitted.
\begin{lem}\label{PSalpha}
\begin{itemize}
\item[(i)] There exists a $(PS)_{\alpha_{\lambda, k}} $ sequence in $ \mathcal{M}_{\lambda, k} $ for $ J_{\lambda, k} $;
\item[(ii)] there exists a $(PS)_{\alpha^{+}_{\lambda, k}} $ sequence in $ \mathcal{M}^{+}_{\lambda, k} $ for $ J_{\lambda, k} $;
\item[(iii)] there exists a $(PS)_{\alpha^{-}_{\lambda, k}} $ sequence in $ \mathcal{M}^{-}_{\lambda, k} $ for $ J_{\lambda, k} $.
\end{itemize}
\end{lem}

\section{ Existence of a ground state solution}

The first lemma in this section establishes the interval where the functional $J_{\lambda, k} $ satisfies the Palais-Smale condition and its statement is the following:

\begin{lem}\label{PS-cond}
Under the assumptions \ref{H1} and \ref{H2}, if $ 0 < \lambda < \Lambda_{1} $, then functional $ J_{\lambda, k} $ satisfies the $(PS)_{d}$ condition for
$$ d < c_{f_{\infty}} - M\left( \lambda^{\frac{p_{+}}{p_{+} - q_{-}}} + \lambda^{\frac{p_{-}}{p_{-} - q_{+}}}  \right) .$$
\end{lem}

\begin{pf}
Let $ \{v_{n}\} \subset W^{1,p(x)}(\mathbb{R}^{N}) $ be a $(PS)_{d}$ sequence for functional $ J_{\lambda, k} $ with $ d < c_{f_{\infty}} - M\left( \lambda^{\frac{p_{+}}{p_{+} - q_{-}}} + \lambda^{\frac{p_{-}}{p_{-} - q_{+}}}  \right) $. By Lemma~\ref{ltdalem}, $ \{v_{n}\} $ is a bounded sequence in $ W^{1,p(x)}(\mathbb{R}^{N}) $, and so, for some subsequence, still denoted by $ \{v_{n}\} $,
\[
v_{n} \rightharpoonup v \mbox{ in } W^{1,p(x)}(\mathbb{R}^{N}),
\]
for some $v \in W^{1,p(x)}(\mathbb{R}^{N}).$
Since $ J'_{\lambda, k}(v) = 0 $ and  $ J_{\lambda, k}(v) \geq 0 $, from (\ref{J-J0-on1})-(\ref{J'-J0'-on1}),  $w_{n} = v_{n} - v$ is a $(PS)_{d^{*}}$ sequence for $J_{0, k}$ with
$$
d^*=d - J_{\lambda, k}(v) < c_{f_{\infty}}.
$$

\medskip
\begin{claim} \label{C2} There is $ R > 0 $ such that
\[
\limsup_{n \to \infty} \sup_{y \in \mathbb{R}^{N}} \int_{B_{R}(y)}|w_{n}|^{p(x)} = 0.
\]
\end{claim}
Assuming by a moment the claim, we have
\[
\int |w_{n}|^{r(x)} \to 0.
\]
On the other hand, by (\ref{J'-J0'-on1}),  we know that $ J'_{0, k}(w_{n}) = o_{n}(1) $, then

\begin{align*}
\int \left( | \nabla w_{n}|^{p(x)} + |w_{n}|^{p(x)}\right) = o_{n}(1),
\end{align*}
showing that  $ w_{n} \to 0 $ in $ W^{1,p(x)}(\mathbb{R}^{N}) $.

\bigskip

\noindent\textbf{Proof of Claim \ref{C2}:}
If the claim is not true, for each $ R > 0 $ given, we find $ \eta > 0 $ and $ \{y_{n}\} \subset \mathbb{Z}^{N} $  verifying
\[
\limsup_{n \to \infty}\int_{B_{R}(y_{n})} |w_{n}|^{p(x)} \geq \eta > 0.
\]
Once $w_n \rightharpoonup 0$ in $W^{1,p(x)}(\mathbb{R}^{N}) $, it follows that $\{y_n\}$ is an unbounded sequence.  Setting
\[
\tilde{w}_{n} = w_{n}(\cdot + y_{n}),
\]
we have that $ \{\tilde{w}_{n}\} $ is also a $(PS)_{d^*}$ sequence for $J_{0, k}$, and so, it must be bounded. Then, there are $ \tilde{w} \in W^{1,p(x)}(\mathbb{R}^{N})  $ and a subsequence of $ \{\tilde{w}_{n}\} $, still denoted by itself,  such that
$$
\tilde{w}_{n} \rightharpoonup \tilde{w} \in W^{1,p(x)}(\mathbb{R}^{N})\setminus\{0\}.
$$
Moreover, since $ J'_{0, k}(w_{n}) \phi( \cdot - y_{n}) = o_{n}(1) $ for each $ \phi \in  W^{1,p(x)}(\mathbb{R}^{N})$, it is possible to prove that $ \nabla \tilde{w}_{n}(x) \to \nabla \tilde{w}(x) $ a.e. in $ \mathbb{R}^{N} $. Therefore,
\begin{align*}
\int \left( |\nabla \tilde{w}|^{p(x) - 2} \nabla \tilde{w} \nabla \phi + | \tilde{w} |^{p(x) - 2} \tilde{w} \phi \right) = \int f_{\infty} | \tilde{w} |^{r(x) - 2} \tilde{w} \phi,
\end{align*}
from where it follows that $ \tilde{w} $ is a weak solution of the Problem~(\ref{Pf00}). Consequently, after some routine calculations,
$$
c_{f_{\infty}}  \leq  J_{f_{\infty}}(\tilde{w}) - \frac{1}{r_{-}} J'_{f_{\infty}}(\tilde{w})\tilde{w} \leq \liminf_{n \to \infty} \left\{J_{0, k}(w_{n}) - \frac{1}{r_{-}} J'_{0, k}(w_{n})  w_{n}\right\} = d^*
$$
which is a contradiction. Then, the Claim \ref{C2} is true.
\end{pf}

The next theorem shows both the existence of a ground state and that it lies in $\mathcal{M}_{\lambda, k}^{+}$.

\begin{thm}\label{T2}
Assume that \ref{H1} and \ref{H2} hold. Then, there exists $ 0 < \Lambda_{*} < \Lambda_{1}$, such that for $ \lambda \in (0, \Lambda_{*}) $ problem~(\ref{Plkm}) has at least one ground state solution $ u_{0} $. Moreover, we have that $ u_{0} \in \mathcal{M}_{\lambda, k}^{+}$ and
\begin{align}
J_{\lambda, k}(u_{0}) = \alpha_{\lambda, k} = \alpha^{+}_{\lambda, k} \geq - M\left( \lambda^{\frac{p_{+}}{p_{+} - q_{-}}} + \lambda^{\frac{p_{-}}{p_{-} - q_{+}}}  \right).
\end{align}

\end{thm}
\begin{pf}
By Lemma~\ref{PSalpha}~$(i)$, there is a minimizing sequence $ \{ u_{n} \} \subset \mathcal{M}_{\lambda, k} $ for $ J_{\lambda, k} $ such that
$$
J_{\lambda, k}(u_{n}) = \alpha_{\lambda, k} + o_{n}(1) \quad \mbox{and} \quad J'_{\lambda, k}(u_n) = o_{n}(1).
$$
Since $ c_{f_\infty} > 0 $, there is $ 0 < \Lambda_{*} < \Lambda_{1}$ such that
$$
\alpha_{\lambda, k} < 0 < c_{f_\infty} - M\left( \lambda^{\frac{p_{+}}{p_{+} - q_{-}}} + \lambda^{\frac{p_{-}}{p_{-} - q_{+}}}  \right) \quad \mbox{for any} \quad 0 < \lambda < \Lambda_{*}.
$$
By Lemma~\ref{PS-cond}, there is a subsequence of $ \{u_{n} \} $, still denoted by itself,  and $ u_{0} \in W^{1, p(x)}(\mathbb{R}^{N}) $ such that $ u_{n} \to u_{0} $ in $ W^{1, p(x)}(\mathbb{R}^{N}) $.  Thereby, $ u_{0} $ is a solution of (\ref{Plkm}) and $ J_{\lambda, k}(u_{0}) = \alpha_{\lambda, k} $. We assert that $ u_{0} \in \mathcal{M}^{+}_{\lambda, k} $. Otherwise, since $ \mathcal{M}^{0}_{\lambda, k} = \emptyset $ for $ 0 < \lambda < \Lambda_{*}$, we have $ u_{0} \in \mathcal{M}^{-}_{\lambda, k} $. Hence
\begin{equation} \label{Y1}
\int\lambda g_{k}(x)|u_{0}|^{q(x)} > 0.
\end{equation}
Indeed, if $ 0 = \int\lambda g_{k}(x)|u_{0}|^{q(x)}, $ then
$$
0 = \int\lambda g_{k}(x)|u_{n}|^{q(x)} + o_{n}(1)= \int \left( |\nabla u_{n}|^{p(x)} + |u_{n}|^{p(x)} \right) - \int f_{k}(x) |u_{n}|^{r(x)} + o_{n}(1).
$$
Therefrom,
\begin{align*}
\alpha_{\lambda, k} + o_{n}(1) &=  J_{\lambda, k}(u_{n}) \\
& = \int \frac{1}{p(x)} \left( |\nabla u_{n}|^{p(x)} + |u_{n}|^{p(x)} \right) - \lambda\int \frac{g_{k}(x)}{q(x)}|u_{n}|^{q(x)}  - \int \frac{f_{k}(x)}{r(x)}|u_{n}|^{r(x)} \\
& \geq \frac{1}{p_{+}} \int  \left( |\nabla u_{n}|^{p(x)} + |u_{n}|^{p(x)} \right) - \frac{\lambda}{q_{-}} \int g_{k}(x) |u_{n}|^{r(x)} - \frac{1}{r_{-}} \int f_{k}(x) |u_{n}|^{r(x)} \\
& = \left( \frac{1}{p_{+}} - \frac{1}{r_{-}} \right) \int \left( |\nabla u_{n}|^{p(x)} + |u_{n}|^{p(x)} \right) + o_{n}(1)
\end{align*}
leading to
$$
\alpha_{\lambda, k} \geq \left( \frac{1}{p_{+}} - \frac{1}{r_{-}} \right) \limsup_{n \in \mathbb{N}}\int \left( |\nabla u_{n}|^{p(x)} + |u_{n}|^{p(x)} \right)
$$
which is absurd, because $ \alpha_{\lambda, k} < 0 $, showing that (\ref{Y1}) holds.

By Lemma~\ref{Max-Nehari}~(ii), there are numbers $ t^{+} < t^{-} = 1 $ such that $ t^{+}u_{0} \in \mathcal{M}^{+}_{\lambda, k} $, $ t^{-}u_{0} \in \mathcal{M}^{-}_{\lambda, k} $ and
\begin{align*}
J_{\lambda, k} (t^{+} u_{0}) < J_{\lambda, k} (t^{-} u_{0}) = J_{\lambda, k} (u_{0}) = \alpha_{\lambda, k} ,
\end{align*}
which is a contradiction. Thereby, $ u_{0} \in \mathcal{M}^{+}_{\lambda, k} $ and
\begin{align*}
-M \left( \lambda^{\frac{p_{+}}{p_{+} - q_{-}}} + \lambda^{\frac{p_{-}}{p_{-} - q_{+}}}   \right) \leq J_{\lambda, k} (u_{0}) = \alpha_{\lambda, k} = \alpha^{+}_{\lambda, k}.
\end{align*}
\end{pf}

\section{Existence of $ \ell $ solutions}
In this section, we will show that \eqref{Plkm} has at least $ \ell $ nontrivial solutions belonging to $ \mathcal{M}^{-}_{\lambda, k} $.
\subsection{Estimates involving the minimax levels}

The main goal of this subsection is to prove some estimates involving the minimax levels $c_{\lambda, k},c_{0, k}, c_{\infty} $ and
$c_{f_{\infty}}$.

First of all, we recall the inequalities
$$
J_{\lambda, k} (u) \leq J_{0, k}(u) \,\,\, \mbox{and} \,\,\, J_{\infty}(u) \leq J_{0, k}(u) \,\,\,\,\, \forall u \in W^{1,p(x)}(\mathbb{R}^{N}),
$$
which imply
\[
c_{\lambda, k} \leq c_{0, k}\quad  \mbox{ and }  \quad c_{\infty} \leq c_{0, k}.
\]

\begin{lem}\label{c0<cf00}
The minimax levels $c_{0, k}$ and $c_{f_{\infty}}$ satisfy the inequality  $c_{0, k} < c_{f_{\infty}}$. Hence, $c_{\infty} < c_{f_{\infty}}$.
\end{lem}
\noindent \textbf{Proof.}
In a manner analogous to Theorem~\ref{TeoComp}, there is $ U \in W^{1,p(x)}(\mathbb{R}^{N}) $ verifying
\[
J_{f_{\infty}}(U) = c_{f_{\infty}} \quad \mbox{ and } \quad J'_{f_{\infty}}(U) = 0.
\]
Similar to Lemma~\ref{Max-Nehari}, there exists $ t > 0 $ such that $ t U \in \mathcal{M}_{0, k} $. Thus,
$$
c_{0, k} \leq  J_{\infty}(tU) = \int \frac{t^{p(x)}}{p(x)} \left( |\nabla U|^{p(x)} + |U|^{p(x)}\right) - \int  \frac{t^{r(x)}}{r(x)}  |U|^{r(x)}.
$$
By \ref{H2}, $f_{\infty} < f(x) $ for all $ x \in \mathbb{R}^{N}$, and so, $f_\infty <1$. Then,  
$$
c_{0, k} <  J_{f_{\infty}}(tU) \leq \max_{s \geq 0}J_{f_{\infty}}(sU) = J_{f_{\infty}}(U) = c_{f_{\infty}}.
$$
\fim


In what follows, let us fix $ \rho_{0}, r_{0} > 0 $ satisfying
\begin{itemize}
  \item $ \overline{B_{\rho_{0}}(a_{i})} \cap \overline{B_{\rho_{0}}(a_{j})} = \emptyset $ for $ i \neq j$ \,\,\, \mbox{and} \,\,\, $i,j \in \{1,...,\ell\}$
  \item $ \bigcup^{\ell}_{i = 1}B_{\rho_{0}}(a_{i}) \subset B_{r_0}(0) $.

  \item $K_{\frac{\rho_{0}}{2}} = \bigcup^{\ell}_{i = 1}\overline{B_{\frac{\rho_{0}}{2}}(a_{i})}$
\end{itemize}
Furthermore,  we define the function $ Q_{k} : W^{1, p(x)}(\mathbb{R}^{N}) \to \mathbb{R}^N $ by
\begin{align*}
Q_{k}(u) = \frac{\int \chi (k^{-1} x)|u|^{p_{+}}}{\int |u|^{p_{+}}},
\end{align*}
where  $ \chi : \mathbb{R}^{N}  \to \mathbb{R}^{N} $ is given by
\[ \chi(x) = x  \,\mbox{ if  } \,|x| \leq r_{0}  \,\mbox{ and } \,  \chi(x) = r_{0} \frac{x}{|x|}  \, \mbox{ if }\,\,  |x| > r_{0} .
\]

The next two lemmas give important information on the function $Q_{k}$ and the level $c_{\infty}$. 

\begin{lem}\label{lemK}
There are $ \delta_{0} > 0 $ and $ k_1 \in \mathbb{N} $ such that if $ u \in \mathcal{M}_{0, k} $ and $ J_{0, k}(u) \leq c_{\infty} + \delta_{0} $, then
\[
Q_{k}(u) \in K_{\frac{\rho_{0}}{2}} \,\,\,\, \mbox{for} \,\,\, k \geq  k_1.
\]
\end{lem}
\noindent \textbf{Proof.}
If the lemma does not occur, there must be $ \delta_{n} \to 0 $, $ k_{n} \to +\infty $ and $ u_{n} \in \mathcal{M}_{0, k_{n}} $ satisfying
\[
J_{0, k_{n}}(u_{n}) \leq c_{\infty} + \delta_{n}
\]
and
\[
Q_{k_{n}}(u_{n}) \not\in K_{\frac{\rho_{0}}{2}}.
\]
Fixing $ s_{n} > 0 $ such that $ s_{n} u_{n} \in \mathcal{M}_{\infty} $, we have that
\[
c_{\infty} \leq J_{\infty}(s_{n} u_{n}) \leq J_{0, k_{n}} (s_{n}u_{n}) \leq \max_{t \geq 0 } J_{0, k_{n}} (tu_{n}) = J_{0, k_{n}}(u_{n}) \leq c_{\infty} + \delta_{n}.
\]
Hence,
\[\{s_{n} u_{n}\} \subset \mathcal{M}_{\infty} \,\,\,\, \mbox{and} \,\,\,\, J_{\infty}(s_{n} u_{n}) \to c_{\infty}.
\]

Applying the Ekeland's variational principle, we can assume without loss of generality that
$\{s_{n} u_{n}\} \subset \mathcal{M}_{\infty} $  is a  $ (PS)_{c_{\infty}} $  sequence  for $ J_{\infty} $, that is,
$$
J_\infty(s_n u_n) \to c_\infty \,\,\,\, \mbox{and} \,\,\,\, J'_{\infty}(s_n u_n) \to 0.
$$
From Theorem \ref{TeoComp}, we must consider the ensuing cases:
\begin{description}
  \item[i)] $ s_{n}u_{n} \to U \neq 0 $ in $ W^{1,,p(x)}(\mathbb{R}^{N}) $; \par
\end{description}
or
\begin{description}
\item[ii)] There exists $ \{y_{n}\} \subset \mathbb{Z}^{N} $ with $|y_n| \to +\infty $ such that $ v_{n}(x) = s_{n}u(x + y_{n}) $ is convergent in $ W^{1,,p(x)}(\mathbb{R}^{N}) $ for some $ V \in  W^{1,p(x)}(\mathbb{R}^{N}) \setminus \{0\}$.
\end{description}

By a direct computation, we can suppose that $ s_{n} \to s_{0} $ for some $ s_{0} > 0 $. Therefore, without loss of generality,
we can assume that
$$
u_{n} \to U  \,\,\, \mbox{or} \,\,\,\, v_{n} = u( \,\, \cdot + y_{n}) \to  V \,\,\,\, \mbox{in} \,\,\,  W^{1, p(x)}(\mathbb{R}^{N}).
$$

\bigskip

\noindent\textbf{Analysis of} $\mathbf{i)}$.

\bigskip

By Lebesgue's dominated convergence theorem
\[
Q_{k_{n}}(u_{n}) = \frac{\int \chi({k_{n}}^{-1} x)|u_{n}|^{p_{+}}}{\int |u_{n}|^{p_{+}}} \to \frac{\int \chi(0)|U|^{p_{+}}}{\int |U|^{p_{+}}} = 0,
\]
implying $ Q_{k_{n}}(u_{n}) \in K_{\frac{\rho_{0}}{2}} $ for $n $ large, because $0 \in K_{\frac{\rho_{0}}{2}}$. However, this a contradiction, because  we are supposing  $ Q_{k_{n}}(u_{n}) \notin K_{\frac{\rho_{0}}{2}} $ for all $n$.

\bigskip

\noindent\textbf{Analysis of} $\mathbf{ii)}$.

\bigskip
Using again the Ekeland's variational principle, we can suppose that $ J'_{0, k_{n}}(u_{n}) = o_{n}(1) $. Hence, $ J'_{0, k_{n}} (u_{n}) \phi(\cdot - y_{n}) = o_{n}(1) $ for  any $ \phi \in W^{1,,p(x)}(\mathbb{R}^{N})$, and so,
\begin{equation}
o_{n}(1)  = \int \left(  |\nabla v_{n}|^{p(x) - 2} \nabla v_{n} \nabla \phi + |v_{n}|^{p(x) - 2} v_{n} \phi\right) - \int f({k_{n}}^{-1} (x + y_{n}))  |v_{n}|^{r(x) - 2}v_{n} \phi. \label{eq5}
\end{equation}
The last limit implies that for some subsequence,
$$
\nabla v_{n}(x) \to \nabla V(x) \,\,\, \mbox{and} \,\,\,  v_{n}(x) \to  V(x) \,\,\, \mbox{a.e in} \,\,\, \mathbb{R}^{N}.
$$
Now, we will study two cases:
\begin{description}
  \item[I)] $ |{k_{n}}^{-1}y_{n}| \to +\infty $
\end{description}
and
\begin{description}
  \item[II)] $ {k_{n}}^{-1}y_{n} \to y $, for some $y \in \mathbb{R}^{N}$.
\end{description}

If I) holds, we see that
\[
\int \left(  |\nabla V|^{p(x) - 2} \nabla V \nabla \phi + |V|^{p(x) - 2} V \phi\right) = \int f_{\infty} |V|^{r(x) - 2}  V \phi,
\]
showing that $ V $ is a nontrivial weak solution of the problem~(\ref{Pf00}). Now, combining the condition $ f_{\infty} < 1 $  with Fatou's Lemma, we get
$$
c_{f_{\infty}} \leq J_{f_{\infty}}(V)  = J_{f_{\infty}}(V) - \frac{1}{r_{-}}J'_{f_{\infty}}(V)V \leq \liminf_{n \to \infty} \left\{  J_{\infty}(u_{n}) - \frac{1}{r_{-}}J'_{\infty}(u_{n})u_{n}  \right\} = c_{\infty},
$$
or equivalently, $ c_{f_{\infty}} \leq c_{\infty} $, contradicting the Lemma~\ref{c0<cf00}.

\bigskip



Now, if  $ {k_{n}}^{-1}y_{n} \to y $ for some $y \in \mathbb{R}^{N}$, then $ V $ is a weak solution of the following problem
\begin{align}
\left\{
\begin{array}
[c]{rcl}%
-\Delta_{p(x)} u + |u|^{p(x) - 2} u & = & f(y)|u|^{r(x) - 2}u  , \quad
\mathbb{R}^{N}\\
\mbox{}\\
u \in W^{1, p(x)}(\mathbb{R}^{N}). &  &
\end{array}
\right. \tag{$ P_{f(y)} $}\label{Pfy}%
\end{align}
Repeating the previous arguments, we deduce that
\begin{align}
c_{f(y)} \leq c_{\infty}, \label{eq6}
\end{align}
where $ c_{f(y)} $ the mountain pass level of the functional  $ J_{f(y)} : W^{1,p(x)}(\mathbb{R}^{N}) \to \mathbb{R} $ given by
\begin{align*}
J_{f(y)}(u)   =\int \frac{1}{p(x)}\left(  |\nabla
u|^{p(x)}+|u|^{p(x)}\right)    -  \int \frac{ f(y) }{r(x)} |u|^{r(x)}.
\end{align*}
If $ f(y) < 1 $, a similar argument explored in the proof of Lemma~\ref{c0<cf00} shows that $ c_{f(y)} > c_{\infty} $, contradicting the inequality (\ref{eq6}). Thereby, $ f(y) = 1$ and $ y = a_{i} $ for some $ i = 1, \cdots \ell $.
Hence,
\begin{align*}
Q_{k_{n}}(u_{n})  &=  \frac{\int \chi({k_{n}}^{-1} x)|u_{n}|^{p_{+}}}{\int |u_{n}|^{p_{+}}}\\
 & =  \frac{\int \chi({k_{n}}^{-1} x + {k_{n}}^{-1} y_{n})|v_{n}|^{p_{+}}}{\int |v_{n}|^{p_{+}}} \to  \frac{\int \chi(y)|V|^{p_{+}}}{\int |V|^{p_{+}}}=a_i, \\
\end{align*}
implying that $ Q_{k_{n}}(u_{n}) \in K_{\frac{\rho_{0}}{2}} $ for $ n $ large, which is a contradiction, since by assumption $ Q_{k_{n}}(u_{n}) \not\in K_{\frac{\rho_{0}}{2}} $.
\fim

\begin{lem}\label{lemK2}
Let $ \delta_{0}>0 $ given in Lemma \ref{lemK} and $ k_3=\max\{k_1,k_2\}$. Then, there is $ \Lambda^{*} = \Lambda^{*}(k) > 0 $ such that
\[
Q_{k}(u) \in K_{\frac{\rho_{0}}{2}}, \,\,\,\,\,\, \forall (u, \lambda, k) \in \mathcal{A}_{\lambda, k} \times [0, \Lambda^*) \times ( [k_3,+\infty) \cap \mathbb{N}),
\]
where $ \mathcal{A}_{\lambda, k} := \left\{ u \in \mathcal{M}^{-}_{\lambda, k}\colon J_{\lambda, k}(u) < c_{\infty} + \frac{\delta_{0}}{2} \right\} $.
\end{lem}

\noindent \textbf{Proof.}
Observe that
$$
J_{\lambda, k}(u) = J_{0, k}(u) - \lambda \int \frac{g_{k}(x)}{q(x)} |u|^{q(x)} \,\,\, \forall u \in W^{1,p(x)}(\mathbb{R}^{N}).
$$
In what follows, let $ t_{u} > 0 $ such that $ t_{u}u \in \mathcal{M}_{0, k} $. Then,
\begin{align}
J_{0, k}(t_{u} u) & = J_{\lambda, k}(t_{u}u) + \lambda \int \frac{g_{k}(x)}{q(x)} (t_{u})^{q(x)} |u|^{q(x)} \nonumber\\
    & \leq  \max_{t \geq 0} J_{\lambda, k}(tu)+ \lambda \int \frac{g_{k}(x)}{q(x)}  (t_{u})^{q(x)} |u|^{q(x)}.\label{Qk}
\end{align}

\begin{claim}\label{Alke} \quad
\medskip

\noindent \textbf{a)} Given $ \Lambda > 0 $, there is a constant $ R > 0 $ such that
$
\mathcal{A}_{\lambda, k}  \subset B_{R}(0),
$
for all $k \geq k_1$ and $ \lambda \in [0,\Lambda] $,  that is, $ \mathcal{A}_{\lambda, k} $ is bounded set, where $k_1$ was given in Lemma \ref{lemK}. Moreover, $ R $ is independent of  $k$.

\bigskip

\noindent \textbf{b)} Let $ u \in \mathcal{A}_{\lambda, k} $ and $ t_{u} > 0 $ such that $ t_{u} u \in \mathcal{M}_{0, k} $. Then, given $\Lambda >0$, there are $ C > 0 $ and $k_2 \in \mathbb{N}$ such that
\[
0 \leq t_{u} \leq C, \quad \mbox{ for all } (u, \lambda, k)  \in \mathcal{A}_{\lambda, k} \times [0,\Lambda] \times ([ k_2, +\infty) \cap \mathbb{N}).
\]
\end{claim}
\noindent \textbf{Proof of a)}
Let $ u \in \mathcal{M}^{-}_{\lambda, k} \subset \mathcal{M}_{\lambda, k} $ such that $ J_{\lambda, k}(u) < c_{\infty} + \frac{\delta_{0}}{2}$ for $ k \geq k_1 $. Then,
\begin{align*}
\int \left( |\nabla u|^{p(x)} + |u|^{p(x)} \right) - \lambda \int g_{k}(x) |u|^{q(x)} - \int f_{k}(x)|u|^{r(x)} = 0
\end{align*}
and
\begin{align*}
\int  \frac{1}{p(x)} \left( |\nabla u|^{p(x)} + |u|^{p(x)} \right) - & \lambda \int  \frac{g_{k}(x)}{q(x)} |u|^{q(x)} - \int   \frac{f_{k}(x)}{r(x)} |u|^{r(x)}  < c_{\infty} + \frac{\delta_{0}}{2}.
\end{align*}
Combining the last two expressions, we obtain
\begin{align*}
\left( \frac{1}{p_{+}} - \frac{1}{r_{-}} \right) \int  \left( |\nabla u|^{p(x)} + |u|^{p(x)} \right) + \left( \frac{1}{q_{-}} - \frac{1}{r_{-}} \right)\lambda \int g_{k}(x) |u|^{q(x)}
 < c_{\infty} + \frac{\delta_{0}}{2}.
\end{align*}
By previous calculations, we have
\begin{align*}
 \left( \frac{1}{p_{+}} - \frac{1}{r_{-}} \right) & \min\{ \Vert u \Vert^{p_{-}}, \Vert u \Vert^{p_{+}} \}  - \Lambda \left( \frac{1}{q_{-}} - \frac{1}{r_{-}} \right) 2 \Vert g_{k} \Vert_{\Theta(x)} K^{q_{+}} \max\{ \Vert u \Vert^{q_{-}}, \Vert u \Vert^{q_{+}} \} \\
&< c_{\infty} + \frac{\delta_{0}}{2}.
\end{align*}
Since $ q_{+} < p_{-} $, it follows that there is $ R > 0 $ such that
\begin{align*}
\Vert u \Vert \leq R \quad \mbox{ for  all } (u, \lambda, k) \in \mathcal{A}_{\lambda, k} \times [0, \Lambda] \times ([ k_1, +\infty) \cap \mathbb{N})
\end{align*}
proving a).

\bigskip

\noindent \textbf{Proof of b)}
Supposing by contradiction that the lemma does not hold. Then, there is $ \{u_{n}\} \subset \mathcal{A}_{\lambda_n, k_n} $ with $\lambda_n \to 0$ and $k_n \to +\infty$ such that $ t_{u_{n}} u_{n} \in \mathcal{M}_{0, k_n} $ and $ t_{u_{n}} \to \infty $ as $ n \to \infty $. Without loss of generality, we can assume that $ t_{u_{n}} \geq 1$. Since $ t_{u_{n}} u_{n} \in \mathcal{M}_{0, k_n} $ and $f_{\infty} < f(x) $ for all $ x \in \mathbb{R}^{N}$, we derive
$$
(t_{u_{n}})^{p_+} \int  \left( |\nabla u_{n}|^{p(x)} + |u_{n}|^{p(x)} \right) \geq f_{\infty} (t_{u_{n}})^{r_{-}}  \int |u_{n}|^{r(x)},
$$
or equivalently,
\begin{align}
\int \left( |\nabla u_{n}|^{p(x)} + |u_{n}|^{p(x)} \right) \geq  f_{\infty} t_{u_{n}}^{r_{-}- p_{+}}  \int |u_{n}|^{r(x)}   \label{modular-tu}
\end{align}
for $ n $ large enough.

Now, we claim that there is $ \eta_{1} > 0 $ such that
\begin{align}\label{eta1}
\int |u_n|^{r(x)} > \eta_{1} \,\,\, \forall n \in \mathbb{N}.
\end{align}

Indeed, arguing by contradiction, there is a subsequence, still denoted by $ \{ u_{n} \} $ such that
\begin{align*}
\int |u_n|^{r(x)} = o_{n}(1) \quad \mbox{ as } n \to \infty.
\end{align*}
As $ u_n \in  \mathcal{M}^{-}_{\lambda_n, k_n}  \subset \mathcal{M}_{\lambda_n, k_n} $, we get
\begin{align*}
(p_{-} - q_{+}) \int \left( |\nabla u_{n}|^{p(x)} + |u_{n}|^{p(x)} \right) - (r_{+} - q_{+}) \int f_{k}(x) |u_{n}|^{r(x)}  < 0.
\end{align*}
By item~a), there are positive constants $ c_{1} $ and $ c_{2} $ such that $ c_{1} < \rho_{1} (u_n) < c_{2} $. Thus,
\begin{align*}
\frac{p_{-} - q_{+}}{r_{+} - q_{+}} < \frac{\int f_{k}(x) |u_n|^{r(x)} }{\int \left( |\nabla u_{n}|^{p(x)} + |u_{n}|^{p(x)} \right)}< \frac{\int |u_{n}|^{r(x)}}{c_{1}} = o_{n}(1)
\end{align*}
which is a contradiction, proving the claim. Thereby, from inequality (\ref{modular-tu}),
$$
\rho_{1}(u_{n})=\int \left( |\nabla u_{n}|^{p(x)} + |u_{n}|^{p(x)} \right)  \to +\infty,
$$
implying that $ \{u_{n}\}$ is a unbounded sequence. However, this is impossible, because by item~a), $ \{u_{n} \} $ is bounded, showing that b) holds.

\medskip

By Claim~\ref{Alke}-b and H\"{o}lder's inequality, It follows of \eqref{Qk} that
\begin{align*}
J_{0, k}(t_{u} u) \leq J_{\lambda, k}(u) + \frac{\lambda}{q_{-}} C^{q_+} \Vert g_{k} \Vert_{\Theta(x)} \Vert |u|^{q(x)}\Vert_{\frac{r(x)}{q(x)}}.
\end{align*}
Once that $u \in \mathcal{A}_{\lambda, k}$, we get
\begin{align*}
J_{0, k}(t_{u} u) < c_{\infty} + \frac{\delta_{0}}{2} + \lambda c_{2} \Vert g_{k} \Vert_{\Theta(x)}  \Vert |u|^{q(x)}\Vert_{\frac{r(x)}{q(x)}}.
\end{align*}
Using the Sobolev embedding combined with Claim~\ref{Alke}-a), we obtain
\[
J_{0, k} (t_{u} u) < c_{\infty} + \frac{\delta_{0}}{2} +  \lambda c_{3}\Vert g_{k} \Vert_{\Theta(x)} \quad \forall u \in \mathcal{A}_{\lambda, k}
\]
where $c_{3} $ is a positive constant. Setting $ \Lambda^{*} : = {\delta_{0}}/{2c_{3}\Vert g_{k} \Vert_{\Theta(x)}}$  and $\lambda \in [0, \Lambda^*)$, we get
\[
t_{u}u \in \mathcal{M}_{0, k} \quad \mbox{ and } \quad J_{0, k} (t_{u} u) < c_{\infty} + \delta_{0}.
\]
Then, by Lemma~\ref{lemK},
\[
Q_{k}(t_{u} u) \in K_{\frac{\rho_{0}}{2}}.
\]
Now, it remains to note that
\[
Q_{k}(u) = Q_{k}(t_{u} u),
\]
to conclude the proof of lemma. \fim

\vspace{0.5 cm}

From now on, we will use the ensuing notations
\begin{itemize}
  \item $ \theta^{i}_{\lambda, k} = \left\{u \in \mathcal{M}^{-}_{\lambda, k} ; |Q_{k}(u) - a_{i} | < \rho_{0} \right\}$,
  \item $ \partial\theta^{i}_{\lambda, k} = \left\{u \in \mathcal{M}^{-}_{\lambda, k} ; |Q_{k}(u) - a_{i} | = \rho_{0} \right\}$,
  \item $ \beta^{i}_{\lambda, k} = \displaystyle\inf_{u \in \theta^{i}_{\lambda, k}} J_{\lambda, k}(u) $
  \end{itemize}
and
\begin{itemize}
 \item  $ \tilde{\beta}^{i}_{\lambda, k} = \displaystyle\inf_{u \in \partial\theta^{i}_{\lambda, k}} J_{\lambda, k}(u) .$
\end{itemize}

The above numbers are very important in our approach, because we will prove that there is a $(PS)$ sequence of $J_{\lambda, k}$ associated with each $\theta^{i}_{\lambda, k}$ for $i=1,2,...,\ell$. To this end, we need of the following technical result

\begin{lem}\label{lemma_A}
There are $ 0 < \Lambda_{\sharp} < \Lambda^{*} $ and $ k \geq k_{\sharp} $ such that
\begin{align*}
\beta_{\lambda, k} < c_{f_{\infty}} - M \left( \lambda^{\frac{p_{+}}{p_{+} - q_{-}}} + \lambda^{\frac{p_{-}}{p_{-} - q_{+}}}  \right) \quad \mbox{ and } \quad \beta^{i}_{\lambda, k} < \tilde{\beta}^{i}_{\lambda, k}
\end{align*}
for all $ \lambda \in [0, \Lambda_{\sharp}) $ and $ k \geq k_{\sharp} $.
\end{lem}
\begin{pf}
From now on,  $ U \in W^{1, p(x)}(\mathbb{R}^{N}) $ is a ground state solution associated with (\ref{Poo}), that is,
\[
J_{\infty}(U) = c_{\infty} \quad \mbox{ and } \quad  J'_{\infty}(U) = 0 \,\,\, \mbox{( See Theorem \ref{TeoComp} )}.
\]
For $ 1 \leq i \leq \ell$ and $ k \in \mathbb{N} $, we define the function $ \widehat{U}^{i}_{k} : \mathbb{R}^{N} \to \mathbb{R}$ by
\[
\widehat{U}^{i}_{k}(x) = U( x - ka_{i}).
\]

\begin{claim}\label{ltda-Jl-coo} For all $i \in \{1,...,\ell \}$, we have that
\[
\limsup_{k \to +\infty}\Bigl(\sup_{t \geq 0 } J_{\lambda, k}(t\widehat{U}^{i}_{k})\Bigr) \leq c_{\infty}.
\]
\end{claim}
\noindent Indeed, since $p, q$ and $ r $ are $ \mathbb{Z}^{N} $-periodic, and $ a_i \in \mathbb{Z}^{N} $, a making a change variable gives
\begin{align*}
J_{\lambda, k}(t\widehat{U}^{i}_{k}) = &\int  \frac{t^{p(x)}}{p(x)} \left( |\nabla U|^{p(x)} + |U|^{p(x)} \right)
- \lambda  \int g(k^{-1} x + a_{i})\frac{t^{q(x)}}{q(x)} \left| U \right|^{q(x)}\\
& \qquad - \int f(k^{-1} x + a_{i}) \frac{t^{r(x)}}{r(x)} |U|^{r(x)}.
\end{align*}
Moreover, we know that there exists $ s = s(k)  > 0 $ such that
\begin{align*}
\max_{t \geq 0} J_{\lambda, k}(t\widehat{U}^{i}_{k}) = J_{\lambda, k}(s\widehat{U}^{i}_{k}) \geq \beta,
\end{align*}
where $\beta$ was given in Lemma \ref{local}. By a direct computation, it is possible to prove that
$$
s(k) \not\to 0  \quad \mbox{and} \quad  s(k) \not\to \infty  \quad \mbox{as} \quad k \to \infty.
$$
Thus, without loss of generality, we can assume $ s(k) \to s_{0}>0$ as $ k \to \infty $. Thereby,
\begin{align*}
\lim_{k \to \infty} \Bigl( \max_{t \geq 0}J_{\lambda, k}(t\widehat{U}^{i}_{k})  \Bigr) & \leq \int \frac{s_{0}^{p(x)}}{p(x)} \left( |\nabla U|^{p(x)} + |U|^{p(x)} \right) - \int f(a_{i}) \frac{s_{0}^{r(x)}}{r(x)} |U|^{r(x)}\\
& \leq  J_{\infty}(s_{0}U) \leq \max_{s \geq 0}J_{\infty} (sU) = J_{\infty}(U) = c_{\infty}.
\end{align*}
Consequently,
\begin{align*}
\limsup_{k \to +\infty}\big(\sup_{t \geq 0 } J_{\lambda, k}(t\widehat{U})\big) \leq  c_{\infty} \,\,\, \,\, \mbox{for} \,\,\, i \in \{1,....,\ell\},
\end{align*}
showing the claim.

\medskip
By Lemma~\ref{c0<cf00}, there is $ 0 < \Lambda_{\sharp} < \Lambda^{*} $ such that
\begin{align*}
 c_{\infty} < c_{f_{\infty}} - M \left( \lambda^{\frac{p_{+}}{p_{+} - q_{-}}} + \lambda^{\frac{p_{-}}{p_{-} - q_{+}}}  \right) \quad \mbox{ for any } \lambda \in [0, \Lambda_{\sharp}).
\end{align*}
Choosing $ 0 < \bar\delta <  \delta_{0} $ so that
\begin{align*}
c_{\infty} + \bar\delta < c_{f_{\infty}}  - M \left( \lambda^{\frac{p_{+}}{p_{+} - q_{-}}} + \lambda^{\frac{p_{-}}{p_{-} - q_{+}}}  \right) \quad \mbox{ for any } \lambda \in [0, \Lambda_{\sharp}).
\end{align*}
Since $ Q_{k}(U^{i}_{k}) \to a_{i} $ as $ k \to \infty $, then $ U^{i}_{k} \in \theta^{i}_{\lambda, k}  $ for all $ k $ large enough. On the other hand, by Claim~\ref{ltda-Jl-coo}, $ J_{\lambda, k} (U^{i}_{k}) < c_{\infty} + \frac{\bar\delta}{2} $ holds also for $k$ large enough and $\lambda \in [0, \Lambda_\sharp)$. This way, there exists $k_4 \in \mathbb{N}$ such that
$$
\beta^{i}_{\lambda, k} < c_{\infty} + \frac{\bar\delta}{2} < c_{f_{\infty}}  - M \left(\lambda^{\frac{p_{+}}{p_{+} - q_{-}}} + \lambda^{\frac{p_{-}}{p_{-} - q_{+}}}  \right), \,\,\,\, \forall \lambda \in [0, \Lambda_\sharp) \,\,\, \mbox{and} \,\,\, k \geq k_4.
$$

In order to prove the other inequality, we observe that Lemma \ref{lemK2} yields $ J_{\lambda, k}(U^{i}_{k}) \geq c_{\infty} + \frac{\delta_{0}}{2} $ for all $ u \in \partial \theta^{i}_{\lambda, k} $, if $\lambda \in [0, \Lambda_\sharp)$ and $k \geq k_3$. Therefore,
\begin{align*}
\tilde{\beta}^{i}_{\lambda, k} \geq c_{\infty} + \frac{\delta_{0}}{2}, \,\,\, \mbox{for} \,\,\, \lambda \in [0, \Lambda_*) \,\,\,\, \mbox{and} \,\,\, k \geq k_3.
\end{align*}
Fixing $k_\sharp=\max\{k_3,k_4\}$, we derive that
\[
\beta^{i}_{\lambda, k} < \tilde{\beta}^{i}_{\lambda, k},
\]
for $ \lambda \in [0, \Lambda_{\sharp}) $ and $ k \geq k_\sharp$.
\end{pf}

\begin{lem}\label{PSb}
For each $ 1 \leq i \leq\ell$, there exists a $(PS)_{\beta^{i}_{\lambda, k}}$ sequence,   $ \left\{ u^{i}_{n} \right\} \subset \theta^{i}_{\lambda, k} $ for functional $ J_{\lambda, k} $.
\end{lem}
\begin{pf}
By Lemma~\ref{lemma_A}, we know that $ \beta^{i}_{\lambda, k} < \tilde{\beta}^{i}_{\lambda, k} $. Then, the result follows adapting the same ideas explored in \cite{Lin12}.
\end{pf}

\section{Proof of Theorem~\ref{T1} }

Let $ \{ u^{i}_{n} \} \subset \theta^{i}_{\lambda, k} $ be a $(PS)_{\beta^{i}_{\lambda, k}} $ sequence in $ \mathcal{M}^{-}_{\lambda, k} $  for functional $ J_{\lambda, k} $ given by Lemma~\ref{PSb}. Since $ \beta^{i}_{\lambda, k} < c_{f_{\infty}}  - M \left( \lambda^{\frac{p_{+}}{p_{+} - q_{-}}} + \lambda^{\frac{p_{-}}{p_{-} - q_{+}}} \right)$, by Lemma~\ref{PS-cond} there is $ u^{i}$ such that $ u^{i}_{n} \to u^{i} $ in $ W^{1, p(x)}(\mathbb{R}^{N}) $. Thus,
$$
u^{i} \in \theta^{i}_{\lambda, k}, \,\,\, J_{\lambda, k}(u^i) = \beta^{i}_{\lambda} \,\,\, \mbox{and} \,\,\, J'_{\lambda, k}(u^i) = 0.
$$
Now, we infer that $ u^{i} \neq u^{j} $ for $ i \neq j $ as $ 1 \leq i,j \leq \ell $. To see why, it remains to observe that
$$
Q_k(u^i) \in \overline{B_{\rho_{0}}(a_{i})} \,\,\, \mbox{and} \,\,\,\, Q_k(u^j) \in \overline{B_{\rho_{0}}(a_{j})}.
$$
Since
$$
\overline{B_{\rho_{0}}(a_{i})} \cap \overline{B_{\rho_{0}}(a_{j})} = \emptyset \,\,\, \mbox{for} \,\,\, i \not= j,
$$
it follows that $u^i \not= u^j$ for $i \not= j$. From this, $ J_{\lambda, k} $ has at least $ \ell $ critical points in $ \mathcal{M}^{-}_{\lambda, k} $  for $\lambda \in [0, \Lambda_\sharp)$ and $ k \geq k_\sharp$.  By Theorem~\ref{T2} it follows that the problem (\ref{Plkm}) admits at least $ \ell + 1 $ solutions for $\lambda \in [0, \Lambda_\sharp)$ and $ k \geq k_\sharp$.
\fim

\begin{flushright}
\scriptsize{ Claudianor O. Alves \\
 Jos\'e L.P. Barreiro   }\\
\smallskip
  \scriptsize{Universidade Federal de Campina Grande\\
   Unidade Acad\^emica de Matem\'atica\\
   58109-970 Campina Grande, PB - Brasil}

{ \scriptsize emails: coalves@yahoo.com.br \\ lindomberg@dme.ufcg.edu.br}
\end{flushright}

\begin{flushright}
\scriptsize{Jos\'e V.A. Gon\c calves}\\
\smallskip
  \scriptsize{Universidade Federal de Goi\'as\\
   Instituto de Matem\'atica e Estat\'istica\\
   74001-970 Goi\^ania, GO - Brasil}\\

\scriptsize{ email:  goncalves.jva@gmail.com}\\
\end{flushright}

\end{document}